\documentclass[leqno,11pt]{article}
\newcommand{\ve}{\varepsilon}

\newcommand{\qed}{q.e.d.}
\usepackage{latexsym}
\usepackage{amsmath}

\begin{document}

\begin{center}
\begin{Large}\begin{bf}
Three subjects of the Jacobi-Perron algorithm of dimension 2
\footnote{2020 Mathematics Subject Classification ; Primary 11J70 ; Secondary 11A55, 11K50.}  

\footnote{Keywords and phrases ; The Jacobi-Perron algorithm,\ Generalizations of continued fractions,
\ Ideal convergence,\ Approximations by rationals. }

\end{bf}\end{Large}
\end{center}

\bigskip

\begin{center}
{\sc Tsutomu Shimada} 
\end{center}

\bigskip
{\bf 1.}\ Introduction.                     

\medskip

The Jacobi-Perron Algorithm, defined by C.G.J.Jacobi[4] and O.Perron[7], is one 
of the generalizations of the continued fraction algorithm to higher dimensions. 
In this paper we shall treat the case of dimension 2. 
Let $[ \alpha ]$ be the greatest rational integer contained in a real number $\alpha$ 
and $\{ \alpha \}=\alpha-[ \alpha ]$, the fractional part of $\alpha$. 
We denote by $\mathcal{D}$ the domain $\{(x,y)\ ;\ 0 \leq x \leq y,\ 1 \leq y \}$  
and   $\mathcal{D}^{\prime}$ the domain $\{(x,y)\ ;\ 1 \leq x \leq y,\ 1 \leq y \}$.   
Let 
 $\mathcal{D}^{\circ} $ and $(\mathcal{D}^{\prime})^{\circ} $ be the set of all interior points 
 of  $\mathcal{D}$ and $\mathcal{D}^{\prime}$, respectively. 
For $(\alpha_0,\, \beta_0) \in \mathcal{D}$   
 we define $\alpha_n$ and $\beta_n\ (n \geq 1)$ 
to have equalities, 
$\alpha_{n-1}=[\alpha_{n-1}]+\displaystyle{ 1 \over \beta_n}$ and 
$\beta_{n-1}=[\beta_{n-1}]+\displaystyle{ \alpha_n \over \beta_n}\ \ (n \geq 1)$.  
Then $( \alpha_n,\,\beta_n) \in \mathcal{D}$,  
$\alpha_n=\displaystyle{ \{\beta_{n-1} \} \over \{ \alpha_{n-1} \} }$ and 
$\beta_n=\displaystyle{ 1 \over \{ \alpha_{n-1} \} }\ \ (n \geq 1)$.  
We write 
$\psi(\alpha_{n-1},\,\beta_{n-1})=(\alpha_{n},\,\beta_{n})\ (n \geq 1)$ and 
$\psi^n(\alpha_0,\,\beta_0)=(\alpha_{n},\,\beta_{n})\ (n \geq 1)$. 
We call $\psi$ the Jacobi-Perron map. 
Let $a_n=[ \alpha_n ]$,  $b_n=[ \beta_n ]\ \ (n \geq 0)$ and  write  
\[
\left( \!
\begin{array}{@{\,}cc@{\,}}
\alpha_0  \\
\beta_0  
\end{array}
\!\! \right) 
=
\left[
\begin{array}{@{\,}ccccc@{\,}}
a_0 & a_1 & \cdots & a_n & \cdots  \\
b_0 & b_1 & \cdots & b_n & \cdots 
\end{array}
\right]  
\]
the Jacobi-Perron expansion of a point $(\alpha_0,\, \beta_0)$.  
We call $(a_0,\,b_0)$ the first step and $(a_1,\,b_1)$ the second step, and so on. 
When $\{ \alpha_m \}=0$, the expansion stops at $(m+1)$th step.   
It is well-known (Perron [7]) that if $1,\,\alpha_0,\,\beta_0$ are linearly independent 
over ${\bf Q}$, the field of rationals, then 
$\{ \alpha_n \} \not=0$ and $\{ \beta_n \} \not=0$ for all $n \geq 0$ 
and the expansion continues. 
Assume $ \{ \alpha_n \}>0$ and $a_n=b_n$ for some $n \geq 1$, then 
$ \{ \alpha_n \} < \{ \beta_n \}$, $\alpha_{n+1}>1$ and $a_{n+1} \geq 1$ 
because $\alpha_n < \beta_n$. 
Let ${\bf Z}$ be the ring of rational integers. 
A sequence 
$\{(a_n,\,b_n)\,;\,n \geq 0,\,a_n,b_n \in {\bf Z},\,0 \leq a_n \leq b_n,\,1 \leq b_n \}$
is called admissible when 
``$a_{n+1} \geq 1$ if $a_n=b_n\ (n \geq 0)$'' .
Following Bernstein[1], we define  
\[K_n=\left(          
   \begin{array}{ccc}
     0  &   0  &  1   \\
     1  &   0  &  a_n  \\
     0  &   1  &  b_n   
    \end{array}
   \right).  \]
We call $K_n$ a matrix of Jacobi-Perron type. 
We denote by $L_n\,(n \geq 0)$ a product of $n+1$ matrices of Jacobi-Perron type. 
That is $L_n=K_0 K_1 \cdots K_n$ with some $K_i\ (0 \leq i \leq n)$.  
Note that $\det K_n=\det L_n=1$ and  

\medskip                        
$
K_n \!\! \left( \!\!\!         
   \begin{array}{c}
      1   \\
      \alpha_{n+1}  \\
      \beta_{n+1}   
    \end{array}
  \!\!\! \right)
   \!=\beta_{n+1} \!\! \left(  \!\!\!        
   \begin{array}{c}
      1   \\
      \alpha_n  \\
      \beta_n   
    \end{array}
  \!\!\! \right)  \!\!(n \geq 0),\ L_{n-1} \!\!
\left(  \!\!\!        
   \begin{array}{c}
      1   \\
      \alpha_n  \\
      \beta_n   
    \end{array}
  \!\!\! \right)
   \!=\beta_1 \beta_2 \cdots \beta_n \!\!
\left( \!\!\!         
   \begin{array}{c}
      1   \\
      \alpha_0  \\
      \beta_0   
    \end{array}
   \!\!\!\right) \!\! (n \geq 1).  \hfill (1.1)
$  
\medskip

\noindent
We write 
\[ L_n=\left(          
   \begin{array}{ccc}
     r_{n-2}  &   r_{n-1}  &  r_n   \\
     p_{n-2}  &   p_{n-1}  &  p_n   \\
     q_{n-2}  &   q_{n-1}  &  q_n  
    \end{array}
   \right) \  (n \geq 0).
\]
For example 
\[ L_0=K_0=\left(          
   \begin{array}{ccc}
     0  &   0  &  1   \\
     1  &   0  &  a_0  \\
     0  &   1  &  b_0   
    \end{array}
   \right)
=\left(          
   \begin{array}{ccc}
     r_{-2}  &   r_{-1}  &  r_0   \\
     p_{-2}  &   p_{-1}  &  p_0   \\
     q_{-2}  &   q_{-1}  &  q_0  
    \end{array}
   \right), 
\]

\[ L_1=K_0 K_1=\left(          
   \begin{array}{ccc}
     0  &   1     &  b_1   \\
     0  &   a_0  &  1+a_0 b_1  \\
     1  &   b_0  &  a_1+b_0 b_1   
    \end{array}
   \right)
=\left(          
   \begin{array}{ccc}
     r_{-1}  &   r_{0}  &  r_1   \\
     p_{-1}  &   p_{0}  &  p_1   \\
     q_{-1}  &   q_{0}  &  q_1  
    \end{array}
   \right). 
\]

\noindent
It holds that 

\medskip                                   

$    
\hspace{20mm}\left\{ 
\begin{array}{rcl}
 r_n &  =  &  r_{n-1}b_n + r_{n-2}a_n +r_{n-3}       \\
 p_n &  =  &  p_{n-1}b_n + p_{n-2}a_n +p_{n-3} \ \ \ \ \ (n \geq 1).     \\
 q_n &  =  &  q_{n-1}b_n + q_{n-2}a_n +q_{n-3}     
    \end{array} \right. \hfill (1.2)
$
\medskip

\noindent
Equalities (1.1) lead to  

\medskip                             
$
\hspace{1mm}
\left\{ \!\!\!
   \begin{array}{@{\,}l@{\,}}
       \ \ \beta_1 \beta_2 \cdots \beta_n =   r_{n-1}\beta_n + r_{n-2}\alpha_n +r_{n-3},     \\  [8pt]
      \ \ \alpha_0=\displaystyle{ p_{n-1}\beta_n + p_{n-2}\alpha_n +p_{n-3}  
     \over  r_{n-1}\beta_n + r_{n-2}\alpha_n +r_{n-3} },\  
     \ \   \beta_0=\displaystyle{ q_{n-1}\beta_n + q_{n-2}\alpha_n +q_{n-3}  
     \over  r_{n-1}\beta_n + r_{n-2}\alpha_n +r_{n-3} }  \ \  (n \geq 1).
\end{array}  \right.  \hfill(1.3)
$

\medskip 

\noindent
From this, about $\displaystyle{ p_{n}  \over  r_{n} }$ and   $\displaystyle{ q_{n}  \over  r_{n} }$ we have 

\medskip 
 
$ 
\left\{ 
   \begin{array}{@{\,}c@{\,}}
\ \  \min \Big\{ 
\displaystyle{ p_{n-3}  \over  r_{n-3} },\,\displaystyle{ p_{n-2}  \over  r_{n-2} },\,\displaystyle{ p_{n-1}  \over  r_{n-1} } 
\Big\} 
< \displaystyle{ p_{n}  \over  r_{n} } < 
\max \Big\{ 
\displaystyle{ p_{n-3}  \over  r_{n-3} },\,\displaystyle{ p_{n-2}  \over  r_{n-2} },\,\displaystyle{ p_{n-1}  \over  r_{n-1} } 
\Big\}   \\    [12pt]  
\ \  \min \Big\{ 
\displaystyle{ q_{n-3}  \over  r_{n-3} },\,\displaystyle{ q_{n-2}  \over  r_{n-2} },\,\displaystyle{ q_{n-1}  \over  r_{n-1} } 
\Big\} 
< \displaystyle{ q_{n}  \over  r_{n} } < 
\max \Big\{ 
\displaystyle{ q_{n-3}  \over  r_{n-3} },\,\displaystyle{ q_{n-2}  \over  r_{n-2} },\,\displaystyle{ q_{n-1}  \over  r_{n-1} } 
\Big\} 
\end{array}  \right.  \hfill(1.4)  
$   

\medskip 

\noindent  
In this inequalities, we may take  $\displaystyle{ p_{n}  \over  r_{n} }$ and $\displaystyle{ q_{n}  \over  r_{n} }$ for 
$\alpha_0$ and $\beta_0$, respectively.  
Note that neither   
$\displaystyle{ p_{n-3}  \over  r_{n-3} } = \displaystyle{ p_{n-2}  \over  r_{n-2} } = \displaystyle{ p_{n-1}  \over  r_{n-1} } $  
nor  
$\displaystyle{ q_{n-3}  \over  r_{n-3} } = \displaystyle{ q_{n-2}  \over  r_{n-2} } = \displaystyle{ q_{n-1}  \over  r_{n-1} } $  
hold because $\det L_{n-1} =1$. 

\medskip

We say that $(\alpha_0,\,\beta_0)$ has the periodic Jacobi-Perron expansion if there exist rational 
integers $u\ (\geq 0)$ and $v\ (\geq 1)$ with 
$a_{\nu}=a_{\nu+v}$, $b_{\nu}=b_{\nu+v}$ for all $\nu \geq u$.  
We take $u$ and $v$ be the minimum.  
The periodicity of expansion is also defined by the existence of 
 $u\ (\geq 0)$ and $v\ (\geq 1)$ with 
$\alpha_{\nu}=\alpha_{\nu+v}$, $\beta_{\nu}=\beta_{\nu+v}$. 
We know that both are equivalent by Bernstein[1](Theorem 2).  
We call the sequence 
$(a_u,\,b_u),\ \cdots,(a_{u+v-1},\,b_{u+v-1})$  
the period of Jacobi-Perron expansion of $(\alpha_0,\,\beta_0)$ and $v$ the length of the period. 
When $u=0$ the expansion is called purely periodic. 
Perron[7] proved \ 
$\displaystyle\lim_{n \to \infty} \displaystyle{ p_n \over r_n}=\alpha_0$\ \ and\ 
$\displaystyle\lim_{n \to \infty} \displaystyle{ q_n \over r_n}=\beta_0$. 

\medskip

We define $\Delta_n=p_n-\alpha_0 r_n$ and  $\Delta_n^{\prime}=q_n-\beta_0 r_n$  $(n \geq 0)$.  
The Jacobi-Perron expansion of $(\alpha_0,\,\beta_0)$ such that 
$\Delta_n \to 0$ and $\Delta_n^{\prime} \to 0\ (n \to \infty)$ 
is said to be ideally convergent (see Chapter 3 of [1]). 

In section 2, discussing the estimation of $\Delta_n$ (and $\Delta_n^{\prime}$) 
under the conditions of the distribution of the number $n$ with 
$\Delta_{n-1} \Delta_{n} > 0$ (and $\Delta_{n-1}^{\prime} \Delta_{n}^{\prime} > 0$) 
or under the boundedness of $\beta_n\ (n \geq 0)$,   
we show  sufficient conditions for the ideal convergence 
(Theorems $2.1$ --- $2.6$).  
We shall prove Theorem 2.3(II), that was first essentially proved by Paley-Ursell[6], 
by means of an explicit coefficient.  

On the other estimations of $\Delta_n$ (and $\Delta_n^{\prime}$) 
we refer to Schweiger[9],[10],  Dubois-Farhane-Paysant-Le Roux[2]   
and Nakaishi[5] etc.

Section 3 is devoted to the studying on the connection with  the classical continued fractions  
(i.e., of dimension 1), where we shall prove, under some conditions, 
that any real conjugate of 
$\beta_n+\displaystyle{r_{n-2}  \over r_{n-1}} \alpha_n\ (>1)$, 
which is not necessarily an integer, 
belongs to the interval $(-1,\,0)$ for all sufficiently large $n$ (Theorems 3.1 and 3.2). 

In the classical theory of the continued fractions,  the set of real numbers with 
bounded partial quotients is known to be null (Theorem 196 in Hardy-Wright[3]).  
Schweiger generalized this theorem(Satz 3 in [8]).  
In section 4,  
we shall show another proof of this by means of 
calculating the measure of elementary quadrangles and triangles explicitly (see (4.1) and (4.2)),  
that is the set of all $(\alpha_0,\,\beta_0)$ for which $\beta_n\ (n \geq 0)$ 
is bounded from above is null (Theorem 4.4). 

\vspace{1cm} 

{\bf 2.}\  Sufficient conditions for the ideal convergence.   

\medskip

All properties about $\Delta_n$ proved in this section are the case with 
$\Delta_n^{\prime}$. 
By means of (1.3), it holds 
$ \beta_n \Delta_{n-1}+\alpha_n \Delta_{n-2}+\Delta_{n-3}=0$, so 

\medskip                              

$ 
\hspace{25mm} \Delta_n+\{ \beta_n \}\Delta_{n-1}+\{ \alpha_n \} \Delta_{n-2}=0\ \ (n \geq 2)  \hfill (2.1)
$
\medskip

We see that all of 
$\Delta_{n-2},\ \Delta_{n-1}$ and $\Delta_{n}$ never have the same sign for any $n \geq 2$. 

The followings are immediate consequences of the definitions  
\medskip                           

\noindent
$
\ \hspace{10mm}\left\{  \begin{array}{l}
\Delta_0=-\{ \alpha_0 \},\  \Delta_1=\{ \alpha_0 \} \{ \beta_1 \},\  
\Delta_2=-\{ \alpha_0 \}(\{ \alpha_2 \}-\{ \beta_1 \} \{ \beta_2 \}),   \\ [7pt]
\Delta_0^{\prime }=-\{ \beta_0 \}, \ \Delta_1^{\prime }=\{ \beta_0 \}\{ \beta_1 \}\!-\!\{ \alpha_1 \}, \\ [7pt]
\Delta_2^{\prime }=\{ \alpha_1 \}\{ \beta_2 \}\!+\!\{ \beta_0 \}(\{ \alpha_2 \}\!-\!\{ \beta_1 \} \{ \beta_2 \}).  \\
\end{array} \right. \hfill (2.2)
$ 

\bigskip

{\bf Lemma 2.1.}\ \ {\it 
$\{ \alpha_{n-1} \} \{ \beta_n \}+|\{ \beta_{n-1} \}\{ \beta_n \}-\{ \alpha_n \}| <1$ for all $n \geq 1$. 
}
\medskip

{\it Proof.}\ \ $\{ \alpha_{n-1} \} \{ \beta_n \}+\{ \beta_{n-1} \}\{ \beta_n \}-\{ \alpha_n \}
 = 
\displaystyle{ \alpha_{n+1}+\alpha_{n}\alpha_{n+1}-\beta_n  \over \beta_{n}\beta_{n+1} } 
$ 

$ = 
\displaystyle{ \alpha_{n+1}+a_{n}\alpha_{n+1}-b_n+ \{ \alpha_n \}\alpha_{n+1}-\{ \beta_n \} 
\over \beta_{n}\beta_{n+1} } 
= 
\displaystyle{ \alpha_{n+1}+a_{n}\alpha_{n+1}-b_n   
\over \beta_{n}\beta_{n+1} } 
$

$<
\displaystyle{ \alpha_{n+1}+b_{n}\beta_{n+1}-b_n   
\over \beta_{n}\beta_{n+1} } 
= 
\displaystyle{ \alpha_{n+1}+b_n (\beta_{n+1}-1)    
\over b_n \beta_{n+1}+ \alpha_{n+1} }
<1 
$, and 

$\{ \alpha_{n-1} \} \{ \beta_n \}-(\{ \beta_{n-1} \}\{ \beta_n \}-\{ \alpha_n \})
 = 
\displaystyle{ \alpha_{n+1}-\alpha_{n}\alpha_{n+1}+\beta_n  \over \beta_{n}\beta_{n+1} } 
$

$
= 
\displaystyle{ \alpha_{n+1}-a_{n}\alpha_{n+1}+b_n   
\over \beta_{n}\beta_{n+1} } 
= 
\displaystyle{ \alpha_{n+1}-a_{n}\alpha_{n+1}+b_n   
\over b_n \beta_{n+1}+ \alpha_{n+1} } 
$

$
\leq 
\displaystyle{ \alpha_{n+1}+b_n   
\over b_n \beta_{n+1}+ \alpha_{n+1} } <1
$. 
The lemma is proved.  \ \hfill  \qed

\bigskip

This lemma shows that $|\Delta^{\prime}_2| <1$. 
From (2.1), $\Delta_{n-1}+\{ \beta_{n-1} \}\Delta_{n-2}+\{ \alpha_{n-1} \} \Delta_{n-3}=0$.
Substituting this into (2.1),\ we have 

\medskip    

$ \hspace{10mm}
\Delta_n+(\{ \alpha_n \}-\{ \beta_{n-1} \}\{ \beta_n \}) \Delta_{n-2}-\{ \alpha_{n-1} \} 
\{ \beta_n \} \Delta_{n-3}=0. \hfill  (2.3)
$

\bigskip

{\bf Lemma 2.2.}\ \ {\it 
If $\Delta_{n-2} \Delta_{n-1}\!<0\!$ and $\Delta_{n-1} \Delta_n\!<\!0$ 
for a natural number  $n \geq 2$,  then 
$| \Delta_n |\!<\! \{ \beta_n \}| \Delta_{n-1} |$. }

\medskip

{\it Proof.}\ \ 
In our case the equality (2.1) means $| \Delta_n |+\{ \alpha_n \}| \Delta_{n-2} |
=\{ \beta_n \}| \Delta_{n-1} |$, and our result follows immediately.  
\hfill \ \qed

\bigskip

For any $n \geq 2$, let $n_{\ast}$ be the largest number, if it exists, such that 
$n_{\ast} \leq n$ and $\Delta_{n_{\ast}-1}\Delta_{n_{\ast}}>0$. 

\bigskip

{\bf Lemma 2.3.}\ \ {\it 
If $n=n_{\ast}\ (n \geq 2)$, then $|\Delta_n| < \{ \alpha_n \} |\Delta_{n-2} |. $
}

\medskip

{\it Proof.}\ \ From (2.1), 
$| \Delta_{n} |+\{ \beta_{n} \}| \Delta_{{n}-1} |
=\{ \alpha_{n} \}| \Delta_{{n}-2} |$. 
The result follows. 
 \hfill \ \qed

\bigskip

{\bf Lemma 2.4.}\ \ {\it If $n=n_{\ast}+1$, then  

{\rm (i)}\ $| \Delta_n | \leq \{ \alpha_{n_{\ast}} \}\{ \beta_{n_{\ast}+1} \}| \Delta_{n_{\ast}-2} |$ if 
$\{ \beta_{n_{\ast}} \} \{ \beta_{n_{\ast}+1} \}-\{ \alpha_{n_{\ast}+1} \} \geq 0$. 

{\rm (ii)}\ $| \Delta_n |< \{ \alpha_{n_{\ast}+1} \}| \Delta_{{n_{\ast}-2}} |$ if 
$\{ \beta_{n_{\ast}} \} \{ \beta_{n_{\ast}+1} \}-\{ \alpha_{n_{\ast}+1} \} \!<\! 0$  
and  
$\{ \alpha_{n_{\ast}} \}\!<\!\{ \beta_{n_{\ast}} \}$. 

{\rm (iii)}\ 
$|\Delta_n|\!<\!(\{ \alpha_{n_{\ast}+1} \}+(1-\{ \alpha_{n_{\ast}+1} \})\{ \alpha_{n_{\ast}} \}\{ \beta_{n_{\ast}+1} \})
 \cdot \max( |\Delta_{{n_{\ast}-2}}|, |\Delta_{{n_{\ast}-1}}|)$ 
if 
$\{ \alpha_{n_{\ast}} \}>\{ \beta_{n_{\ast}} \}$. 

Summing up above, we have   
$|\Delta_n|<\displaystyle{ 1+ \{ \alpha_{n} \}   \over 2}
 \cdot \max( |\Delta_{n_{\ast}-2}|, |\Delta_{n_{\ast}-1}|)$.  }

\medskip

{\it Proof.}\ \ 
Note that, by definition, $\Delta_{n_{\ast}-2} \Delta_{n_{\ast}-1}<0$  and 
$\Delta_{n_{\ast}} \Delta_{n_{\ast}+1}<0$. From (2.1), 

\medskip   
$\hspace{20mm}
\{ \alpha_{n_{\ast}} \}|\Delta_{n_{\ast}-2}|=|\Delta_{n_{\ast}}|+ \{ \beta_{n_{\ast}} \}|\Delta_{n_{\ast}-1}|
>\{ \beta_{n_{\ast}} \}|\Delta_{n_{\ast}-1}|.  \hfill  (2.4)
$
\medskip

\noindent
(i)   Because  
$\Delta_{n_{\ast}-2} \Delta_{n_{\ast}-1} < 0$, 
from (2.3),  

$|\Delta_n| \leq \max(\{ \alpha_{n_{\ast}} \} \{ \beta_{n_{\ast}+1} \} |\Delta_{n_{\ast}-2}|,\  
(\{ \beta_{n_{\ast}} \}\{ \beta_{n_{\ast}+1} \}-\{ \alpha_{n_{\ast}+1} \}) |\Delta_{n_{\ast}-1}|) $. 

\noindent
By $(2.4)$,  $\{ \alpha_{n_{\ast}} \}\{ \beta_{n_{\ast}+1} \}|\Delta_{n_{\ast}-2}|>
\{ \beta_{n_{\ast}} \}\{ \beta_{n_{\ast}+1} \} |\Delta_{n_{\ast}-1}|$ 

$>
(\{ \beta_{n_{\ast}} \}\{ \beta_{n_{\ast}+1} \}-\{ \alpha_{n_{\ast}+1} \})|\Delta_{n_{\ast}-1}|$. 
So, $|\Delta_n| \leq \{ \alpha_{n_{\ast}} \} \{ \beta_{n_{\ast}+1} \} |\Delta_{n_{\ast}-2}|$. 

\noindent
(ii) \ In this case, (2.3) means 

\medskip    

$\hspace{3mm}
|\Delta_n|=\{ \alpha_{n_{\ast}} \}\{ \beta_{n_{\ast}+1} \}|\Delta_{n_{\ast}-2}|
+(\{ \alpha_{n_{\ast}+1} \}-\{ \beta_{n_{\ast}} \} \{ \beta_{n_{\ast}+1} \})
|\Delta_{n_{\ast}-1}|.   \hfill (2.5)
$

\medskip

Because $\{ \alpha_{n_{\ast}} \} \{ \beta_{n_{\ast}+1} \}
+(\{ \alpha_{n_{\ast}+1} \}-\{ \beta_{n_{\ast}} \}\{ \beta_{n_{\ast}+1} \})
$

$=\{ \alpha_{n_{\ast}+1} \}+\{ \beta_{n_{\ast}+1} \}(\{ \alpha_{n_{\ast}} \}
-\{ \beta_{n_{\ast}} \})
$ 
$<\{ \alpha_{n_{\ast}+1} \}$, 

(2.5) implies 
$|\Delta_n|
<\{ \alpha_{n_{\ast}+1} \} \max (|\Delta_{n_{\ast}-2}|,\  |\Delta_{n_{\ast}-1}| )$. 
Our assumption and (2.4) lead to $|\Delta_{n_{\ast}-2}|>|\Delta_{n_{\ast}-1}|$, therefore 
$|\Delta_n|<\{ \alpha_{n_{\ast}+1} \}|\Delta_{n_{\ast}-2}|$. 

\noindent
(iii)\ From the assumption, 
$\{ \alpha_{n_{\ast}+1} \}=\displaystyle{\{ \beta_{n_{\ast}} \}  \over  \{ \alpha_{{n_{\ast}}} \} }$,  
so  

$\{ \alpha_{n_{\ast}+1} \}-\{ \beta_{n_{\ast}} \} \{ \beta_{n_{\ast}+1} \}=
\{ \beta_{n_{\ast}} \} (\displaystyle{ 1  \over  \{ \alpha_{n_{\ast}} \} }-\{ \beta_{n_{\ast}+1} \})>0$, 
thus (2.5) holds in this case. And 
$\{ \alpha_{n_{\ast}} \}\{ \beta_{n_{\ast}+1} \}+(\{ \alpha_{n_{\ast}+1} \}-\{ \beta_{n_{\ast}} \} \{ \beta_{n_{\ast}+1} \})
=\{ \alpha_{n_{\ast}+1} \}+\{ \alpha_{n_{\ast}} \}\{ \beta_{n_{\ast}+1} \}(1-\{ \alpha_{n_{\ast}+1} \})$. 
From (2.5), we get our inequality. 

The last statement is now clear by   
$\{ \alpha_{n} \} \{ \beta_{n+1} \} = \displaystyle{ \{ \beta_{n+1} \} \over  \beta_{n+1}  } 
< \displaystyle{1 \over 2}\ (n \geq 1)$ and 
$\{ \alpha_{n} \} < \displaystyle{ 1+ \{ \alpha_{n} \}   \over 2} $.  
 \ \hfill \ \qed

\bigskip

{\bf Lemma 2.5.}\ \ {\it It holds that }
\[
| \Delta_n |< \max \Big( \displaystyle{1+ \{ \alpha_n \} \over 2},\ \{ \beta_n \} \Big)
 \cdot \max(|\Delta_{n-3}|,\ |\Delta_{n-2}|,\ |\Delta_{n-1}|)\  (n \geq 3).  
\] 

{\it Proof.}\ \ 
One of the following three cases holds. 

\noindent
(a) ; $\Delta_{n-2} \Delta_{n-1}<0$ and $\Delta_{n-1} \Delta_{n}<0$. 
Then, from Lemma 2.2, $|\Delta_{n}|\! <\! \{ \beta_n \}| \Delta_{n-1}|$. 

\noindent
(b) ; $\Delta_{n-1} \Delta_{n} > 0$. Then, $n=n_{ \ast }$ and from Lemma 2.3, 
$|\Delta_{n}| < \{ \alpha_n \}| \Delta_{n-2}|$. 

\noindent
(c) ; $\Delta_{n-2} \Delta_{n-1} >0$. Then, $n=n_{ \ast }+1$ and from Lemma 2.4,  
$$
|\Delta_n|<\displaystyle{ 1+ \{ \alpha_{n} \}   \over 2}
 \cdot \max( |\Delta_{n-3}|,\ |\Delta_{n-2}|).
$$ 
\noindent
To sum up, we have completed the proof. \ \hfill  \qed

\bigskip

The result of this lemma is a refinement of an inequality in the proof of Theorem 3 in [2]. 
From $(2.2)$ and using Lemma 2.5 repeatedly, we get  

$| \Delta_n |< \max (|\Delta_0|,\ |\Delta_1|,\ |\Delta_2|)=\{ \alpha_0 \}\ \ (n \geq 0)$. 

\noindent
On the other hand, 
$|\Delta_0^{\prime}|<1$ and $|\Delta_1^{\prime}|<1$\ from (2.2), and 
$|\Delta_2^{\prime}|<1$ from Lemma 2.1,   
so    
$| \Delta_n^{\prime} |< \max(|\Delta_0^{\prime}|,\ |\Delta_1^{\prime}|,\ |\Delta_2^{\prime}|)
<1\ \ (n \geq 0)$.  
Then, 
$\Bigl| \displaystyle{ p_n  \over r_n}-\alpha_0 \Bigr|< \displaystyle{ \{ \alpha_0 \} \over r_n}$ 
and 
$\Bigl| \displaystyle{ q_n  \over r_n}-\beta_0 \Bigr|< \displaystyle{ 1 \over r_n}\ \ (n \geq 0)$.  
Since $r_n \to \infty\ (n \to \infty)$, 
we obtain another proof of a theorem of Perron ;  

\medskip

{\bf Theorem}(Perron[7]).\ \ {\it 
$
\displaystyle\lim_{n \to \infty} \displaystyle{ p_n \over r_n}=\alpha_0$
\  and\ \ 
$
\displaystyle\lim_{n \to \infty} \displaystyle{ q_n \over r_n}=\beta_0. 
$
}

\bigskip

{\bf Remark} ; In the case of $a_n=0$ and $b_n=1$ for all $n \geq 0$, 
$r_n$ is less than those of any other case, so is $p_n$ and $q_n$. 
In this case, we have 
$r_0=r_1=r_2=1$ and $r_n=r_{n-1}+r_{n-3}\ (n \geq 3)$. 
Then, by the well-known procedure, we have  
$$
r_n=\displaystyle{ \lambda^{n-2} \over  \delta_1-\delta_2  }
\Bigl\{
\displaystyle{ 2-2 \lambda \over 3 \lambda-2} 
\displaystyle{ \delta_1^{n-1}-\delta_2^{n-1} \over  \lambda^{n-1} }
+
\displaystyle{ 2- \lambda \over 3 \lambda-2} 
\displaystyle{ \delta_1^{n}-\delta_2^{n} \over  \lambda^{n} } 
\Bigr\}
+
\displaystyle{  \lambda^{n+1}  \over 3 \lambda-2} \ \ (n \geq 0), 
$$
\noindent
where $\lambda$ is the real root of $X^3-X^2-1=0$ and 
$\delta_i \ (i=1,2)$ are the roots of 
$X^2+(\lambda-1)X+\lambda (\lambda-1)=0$. 
Here, $\lambda=1.465 \cdots, |\delta_i|=0.826 \cdots$ and 
$\displaystyle{  \lambda^{3}  \over 3 \lambda-2}=1.0739 \cdots$.  
From these, we can deduce the followings ; 

(i)\ $\displaystyle\lim_{n \to \infty} \displaystyle{ r_n \over \lambda^{n-2} } 
= \displaystyle{  \lambda^{3}  \over 3 \lambda-2}$,\  \  
(ii)\ $\displaystyle\lim_{n \to \infty} \displaystyle{  r_n  \over r_{n-1}} = \lambda$,\  \ 
(iii)\ $\lambda^{n-2} < r_n < \lambda^{n-1} \ ( n \geq 3)$.     

\noindent
In the general case of $a_n$ and $b_n\ (n \geq 0)$, from (iii), it holds that 
$\lambda^{n-2} <r_n\ (n \geq 3)$. 
Therefore, in general, we have  
$$ \Bigl| \displaystyle{ p_n  \over r_n} - \alpha_0 \Bigr|< \{ \alpha_0 \} \Bigl( \displaystyle{ 1 \over \lambda} \Bigr)^{n-2} \ \ 
{\rm and} \ \ 
 \Bigl| \displaystyle{ q_n  \over r_n} - \beta_0 \Bigr|<  \Bigl( \displaystyle{ 1 \over \lambda} \Bigr)^{n-2} \ \ ( n \geq 3). 
$$

\bigskip

{\bf Lemma 2.6.}\ \ {\it If $n=n_{\ast}+2 \geq 4$, then  

\medskip

$
|\Delta_n|\!\!<\!\!\left\{\!\!\!
\begin{array}{l}
({\rm i})\ \displaystyle{ 1 \over 2}\{ \beta_{n_{\ast}+2} \} |\Delta_{n_{\ast}-2}|\ \  
if\ \{ \beta_{n_{\ast}} \} \{ \beta_{n_{\ast}+1} \}-\{ \alpha_{n_{\ast}+1} \} \geq 0, \\  [12pt]
({\rm ii})\ \displaystyle{ 1 \over 2}|\Delta_{n_{\ast}-2}|\   
if\ \{ \beta_{n_{\ast}} \} \{ \beta_{n_{\ast}+1} \}-\{ \alpha_{n_{\ast}+1} \} < 0 \  
and\  \{ \alpha_{n_{\ast}} \}<\{ \beta_{n_{\ast}} \},  \\ [12pt]
({\rm iii})\ \displaystyle{ 3 \over 4} \max(|\Delta_{n_{\ast}-2}|,\ |\Delta_{n_{\ast}-1}|) \ \  
if \ \{ \alpha_{n_{\ast}} \}>\{ \beta_{n_{\ast}} \}.  
\end{array} \right.
$

To sum up, we have 
$|\Delta_n|<\displaystyle{ 3 \over 4} \max(|\Delta_{n_{\ast}-2}|,\ |\Delta_{n_{\ast}-1}|)$. }

\medskip

{\it Proof.}\ \ 
From Lemma 2.2,  
$|\Delta_n|<\{ \beta_n \}|\Delta_{n-1}|$. 
Applying Lemma 2.4 to $n-1$, note that $(n-1)_{\ast}=n-2=n_{\ast}$, we obtain  

(i)\  $| \Delta_n |< \{ \beta_n \}| \Delta_{n-1} |<
\{ \beta_n \} \{ \alpha_{n_{\ast}} \} \{ \beta_{n_{\ast}+1} \} | \Delta_{n_{\ast}-2} |<
\displaystyle{ 1 \over 2}\{ \beta_{n_{\ast}+2} \}| \Delta_{n_{\ast}-2} |$, 

(ii)\  $| \Delta_n |\!<\! \{ \beta_n \}| \Delta_{n-1} |<
\{ \beta_n \}\{ \alpha_{n_{\ast}+1} \}| \Delta_{n_{\ast}-2} |
=\{ \beta_{n_{\ast}+2} \}\{ \alpha_{n_{\ast}+1} \}| \Delta_{n_{\ast}-2} |$

$
<\displaystyle{ 1 \over 2}| \Delta_{n_{\ast}-2} |$, 

(iii)\ \ 
$\{ \beta_n \}(\{ \alpha_{n_{\ast}+1} \}\!+\!(1-\{ \alpha_{n_{\ast}+1} \})\{ \alpha_{n_{\ast}} \}\{ \beta_{n_{\ast}+1} \})$

$\!=\!\{ \alpha_{n-1} \}\{ \beta_n \}\!+\!\{ \alpha_{n-2} \}\{ \beta_{n-1} \}\{ \beta_n \}
-\{ \alpha_{n-2} \}\{ \alpha_{n-1} \}\{ \beta_{n-1} \}\{ \beta_n \}$

$<1-(1-\{ \alpha_{n-1} \}\{ \beta_n \})(1-\{ \alpha_{n-2} \}\{ \beta_{n-1} \})
<\displaystyle{ 3 \over 4}$. 

Therefore, 
$| \Delta_n |< \{ \beta_n \}| \Delta_{n-1} |
<\displaystyle{ 3 \over 4}\max(|\Delta_{n_{\ast}-2}|,\ |\Delta_{n_{\ast}-1}|) $.  
\hfill  \qed

\bigskip

{\bf Lemma 2.7.} {\it If  $n \geq n_{\ast}+3 \geq 5$, then 

$| \Delta_n |\!< \!\displaystyle{ 3 \over 4}
\{ \beta_n \} \{ \beta_{n-1} \} \!\cdots \!\{ \beta_{n_{\ast}+3} \} 
\max(|\Delta_{n_{\ast}-2}|, |\Delta_{n_{\ast}-1}|)$. }

\medskip

{\it Proof.}\ \ 
It is clear from Lemmas 2.2 and 2.6.    \hfill  \qed

\bigskip

We denote by $n_i\ ( i \geq 1,\ n_1<n_2< \cdots )$ all of indices $n$ with 
$\Delta_{n-1} \Delta_{n}>0$, if it exists. 
It holds $n_{i-1}+2 \leq n_i \ (i \geq 2)$ generally. 

\bigskip

{\bf Lemma 2.8.}\ \ {\it 
It holds that    

\noindent
$({\rm i})\ \max(|\Delta_{n_i-2}|,\ |\Delta_{n_i-1}|)<\max(|\Delta_{n_{i-1}-2}|,\ |\Delta_{n_{i-1}-1}|)$. 

If $n_{i-1}+4 \leq n_i$, then 

\noindent
$({\rm ii})\ \max(|\Delta_{n_i-2}|,\ |\Delta_{n_i-1}|)$

$<
\displaystyle{ 3 \over 4}
\{ \beta_{n_i-2} \} \{ \beta_{n_i-3} \} \cdots \{ \beta_{n_{i-1}+3} \}
 \cdot \max(|\Delta_{n_{i-1}-2}|,\ |\Delta_{n_{i-1}-1}|)$.  }

\medskip

{\it Proof.}\ \ 
Note that $(n_i-2)_{\ast}=(n_i-1)_{\ast}=n_{i-1}$.  
Taking $n_i-2$ and $n_i-1$ for $n$ in Lemmas 2.2, 2.3, 2.4, 2.6 and 2.7 we can prove the result.    \hfill  \qed

\bigskip

{\bf Theorem 2.1.}\ \ {\it
Suppose there exist infinite number of $i$ 
such that $n_{i-1}+4 \leq n_i$. 
Then, 
$\lim_{n \to \infty} \Delta_n=0$. }

\medskip

{\it Proof.}\ \ 
For each $n \geq 3$, there exists $j=j(n) \geq 1$ with $n_j = n_{\ast}$.  
Because $j \to \infty\ (n \to \infty)$, the number of $i\ (i \leq j)$ tends to infinity as 
$n \to \infty$ and the inequality $({\rm ii})$ of Lemma 2.8 holds infinitely many times 
as $n \to \infty$. 
Thus $\lim_{n \to \infty} \Delta_n=0$.   \hfill  \qed

\bigskip

{\bf Lemma 2.9.}\ \ \  
$\displaystyle\prod_{n \geq 0}\{ \alpha_n \}=0$ 

\medskip

{\it Proof.}\ \ Since $\{ \alpha_{n-1} \}=\displaystyle{ 1  \over  \beta_n }$,  
the first equality of (1.3) means 
$$\{ \alpha_0 \} \{ \alpha_1 \} \cdots \{ \alpha_n \}
=\displaystyle{ 1  \over  r_n \beta_{n+1}+r_{n-1} \alpha_{n+1}+r_{n-2} }.$$  
The lemma holds because $r_n \to \infty$ as $n \to \infty$ by (1.2). 
  \hfill  \qed

\bigskip

{\bf Lemma 2.10.}\ \  {\it 
If $\ \displaystyle\prod_{n \geq 0}x_n=0$ with real numbers $x_n\  (n \geq 0,\ 0<x_n<1)$,\ then 
$\displaystyle\prod_{n \geq 0} \displaystyle{1+x_n  \over 2}=0$. }

\medskip

{\it Proof.}\ \ We may suppose $x_n > \displaystyle{1  \over 2}$ 
without loss of generality. 
Because, in general, 
$\bigl(\displaystyle{1+x  \over 2} \bigr)^3 <x$ for any 
$\displaystyle{1  \over 2}<x<1$, we have 
$\Bigl( \displaystyle\prod_{n \geq 0} \displaystyle{1+x_n  \over 2} \Bigr)^3 \leq 
\displaystyle\prod_{n \geq 0}x_n =0$ and our lemma is proved.   \hfill  \qed

\bigskip

{\bf Lemma 2.11.}\  {\it Any infinite sub-product of $\displaystyle\prod_{n \geq 0}\!\{ \alpha_n \}$\! 
is zero if and only if $\displaystyle\liminf_{n \to \infty} \beta_n \!>\!1$.  }

\medskip

{\it Proof.}\ \ If $\ \displaystyle\liminf_{n \to \infty} \beta_n >1$, then for any $\delta$ 
$(0< \delta < \displaystyle\liminf_{n \to \infty} \beta_n -1)$ 
there exists a natural number $N_{\delta}$ such that $\beta_n \geq 1+\delta$ 
for all $n \geq N_{\delta}$. 
$\{ \alpha_n \}=1/\beta_{n+1} \leq 1/(1+\delta)$, 
so any infinite sub-product is zero.  
Let $x_l\,(l \geq 1)$ be real numbers such that 
$\displaystyle\prod_{l \geq 1} x_l \not= 0\,(0<x_l <1)$. 
Assume $\ \displaystyle\liminf_{n \to \infty} \beta_n =1$, 
then $\ \displaystyle\limsup_{n \to \infty} \{ \alpha_n \} =1$ and 
we can take $n_l$ with $x_l < \{ \alpha_{n_l} \} < 1$ for each $l \geq 1$. 
Therefore we have an infinite sub-product such that 
$\displaystyle\prod_{l \geq 1}\{ \alpha_{n_l} \} \not= 0$.  
  \hfill  \qed

\bigskip

{\bf Theorem 2.2.}\ \  {\it Suppose $\ \displaystyle\liminf_{n \to \infty} \beta_n >1$ and there exist 
infinitely many $n\ (\geq 1)$ satisfying $\Delta_{n-1} \Delta_n >0$, then 
$\displaystyle\lim_{n \to \infty} \Delta_n=0$.  }

\medskip

{\it Proof.}\ \ For any $n$, let $i=i(n)$ be the natural number such that $n_i=n_{\ast}$, 
$n_i$ is the same as before. 
By means of Lemmas 2.3, 2.4, 2.6 and 2.7, we have 
\[
| \Delta_n|
< 
\max (\displaystyle{ 3 \over 4},\  \{ \alpha_{n_i} \},\ 
\displaystyle{1+\{ \alpha_{n_i+1} \}  \over 2})
\cdot \max (|\Delta_{n_i-2}|,\ |\Delta_{n_i-1}|) 
\] 
and, since $(n_i-2)_{\ast}=(n_i-1)_{\ast}=n_{i-1}$, it holds 
$
\max (|\Delta_{n_i-2}|,\ |\Delta_{n_i-1}|)$

$<
\max (\displaystyle{ 3 \over 4},\  \{ \alpha_{n_{i-1}} \},\ 
\displaystyle{1+\{ \alpha_{n_{i-1}+1} \}  \over 2})
\cdot \max (|\Delta_{n_{i-1}-2}|,\ |\Delta_{n_{i-1}-1}|). 
$ 

\noindent
Using this repeatedly, we have 
\[
| \Delta_n|
< 
C_i \max (|\Delta_{n_1-2}|,\ |\Delta_{n_1-1}|)
\]
where 
\[
C_i=\prod_{1 \leq j \leq i}
\max (\displaystyle{ 3 \over 4},\  \{ \alpha_{n_j} \},\ 
\displaystyle{1+\{ \alpha_{n_j+1} \}  \over 2}). 
\] 
The assumption means that $i \to \infty\ (n \to \infty)$, 
so we can deduce $C_i \to 0$ by Lemmas 2.9 --- 2.11 and therefore  $\Delta_n \to 0\  
(n \to \infty)$.   \hfill  \qed

\bigskip

{\bf Lemma 2.12.}\ \  {\it Assume $\beta_n\ (n \geq 0)$ is bounded from above and let $M\ (>2)$  be a 
real number such that $\beta_n \leq M$ for all $n \geq 0$.  Then, 

{\rm (i)} $\ \displaystyle\limsup_{n \to \infty} \displaystyle{\alpha_n  \over  \beta_n} <1$\ \ \ 
$(i.e., \ \displaystyle\limsup_{n \to \infty} \{ \beta_n \} <1)$,   

{\rm (ii)} $\displaystyle{ \alpha_n  \over  \beta_n } \geq \displaystyle{ 1  \over  M^2 }$ for all $n \geq 0$ \ \ \ 
$( i.e., \ \{ \beta_n \} \geq \displaystyle{ 1  \over  M^2 })$, 

{\rm (iii)} $\beta_n \geq 1+\displaystyle{ 1  \over  M^2 }$ for all $n \geq 0$ \ \ \ 
$(i.e., \ \{ \alpha_n \} \leq \displaystyle{ M^2  \over  M^2+1 })$. 
}

\medskip 

{\it Proof.}\ \ From the assumption, $\{ \alpha_n \} \geq 1/M$. 
Let $\delta$ be a real number with $0< \delta < 1/M$ and $\ve$ with $0< \ve < {\delta}^2$. 

(i)\ Suppose 
$\ \displaystyle\limsup_{n \to \infty} \displaystyle{\alpha_n  \over  \beta_n} =1$, 
then there exists a sub-sequence 
$(\alpha_{n_j},\,\beta_{n_j})$\ $(j \geq 1)$ such that 
$\alpha_{n_j}/\beta_{n_j} \to 1\ (j \to \infty)$  and   
there is a natural number $J$ with 
$0\!<\! \beta_{n_j}-\alpha_{n_j} \!< \! \varepsilon$ for all $j \geq J$.  
It holds one of the following two cases ;

(A)\ $a_{n_j}=b_{n_j}$ and 
$\{ \beta_{n_j} \}-\varepsilon < \{ \alpha_{n_j} \} < \{ \beta_{n_j} \}$, 

(B)\ $a_{n_j}+1=b_{n_j}$ and 
$1-\varepsilon < \{ \alpha_{n_j} \} < 1< 1+ \{ \beta_{n_j} \} < 1+\varepsilon$. 

\noindent
If (A) holds, dividing the inequality by $\{ \alpha_{n_j} \}$, we have  
$1< \alpha_{n_j+1} < 1+\displaystyle{\varepsilon  \over  \{ \alpha_{n_j} \} } 
< 1+\displaystyle{ \varepsilon  \over  \delta } < 1+ \delta$. 
So, $ \{ \alpha_{n_j+1} \} < \delta$. 
This is a contradiction.  
If (B) holds, then  
$ \{ \alpha_{n_j} \} > \{ \beta_{n_j} \}$ and we have 
$ \alpha_{n_j+1} < 1$, $ a_{n_j+1} = 0$, 
$ \{ \alpha_{n_j+1} \} = \alpha_{n_j+1} = \{ \beta_{n_j} \}/ \{ \alpha_{n_j} \}
< \varepsilon/(1-\varepsilon) < \delta$  
because   $0<\delta<1/2$ and $0<\delta^2<1/4$. 
This is a contradiction and (i) is proved. 

(ii)\ $\displaystyle{ \alpha_n  \over  \beta_n } \geq \displaystyle{  \{ \alpha_n \}  \over  M } 
\geq 1/M^2$. 

(iii)\ $\beta_n=b_n+ \{ \beta_n \} \geq 1+1/M^2$.    \hfill  \qed

\bigskip

{\bf Theorem 2.3.}\ \  {\it {\rm (I)} Suppose 
$\ \displaystyle\limsup_{n \to \infty} \displaystyle{\alpha_n  \over  \beta_n} <1$  
and  $\ \displaystyle\liminf_{n \to \infty} \beta_n >1$.   
Then, there exists a  real number $c\ (0<c<1)$ and a natural number $N$ such that 
$
\ \ \ |\Delta_n|< \{ \alpha_0 \} c^{n-N} \ \ and \ \ |\Delta_n^{\prime}|< c^{n-N}
$
for all $n \geq N+3$. 

{\rm (II)}   If $\beta_n\ (n \geq 0)$ is bounded from above, then there exists a real number $a\ (>0)$ 
such that 
$
\ \ \ {\bigl|} \displaystyle{p_n  \over  r_n}-\alpha_0 {\bigr|} < \{ \alpha_0 \} c^{-N}/r_n^{1+a}
\ \ \ \ and\ \ \ \ 
{\bigl|} \displaystyle{q_n  \over  r_n}-\beta_0 {\bigr|} <  c^{-N}/r_n^{1+a}\ \ \ 
$
for all $n \geq 3$, where  $c$ and $N$ are the same as in {\rm (I)}.
}
    
\medskip

{\it Proof.}\ \ 
Note that $\displaystyle\limsup_{n \to \infty} \displaystyle{\alpha_n  \over  \beta_n}
=\displaystyle\limsup_{n \to \infty} \{ \beta_n \}$  and  
$1/\displaystyle\liminf_{n \to \infty}  \beta_n  $  
$= \displaystyle\limsup_{n \to \infty} \{ \alpha_n \}$. 

(I) \ 
Because 
$$ 
\max \Bigl( 
\displaystyle{ 1+\limsup_{n \to \infty} \{ \alpha_n \}  \over  2 },\ \limsup_{n \to \infty} \{ \beta_n \}  \Bigr)  <1,
$$ 
there exists a natural number $N$ such that 
$\max \!\Big( \displaystyle{1+ \{ \alpha_n \} \over 2},\ \{ \beta_n \} \Big)<1$ for all $ n \geq N$. 
Lemma 2.5 implies the existence of a real number $\ve\ ( 0< \ve <1)$ such that 
\[
|\Delta_n| < \varepsilon \max(|\Delta_{n-3}|,|\Delta_{n-2}|,|\Delta_{n-1}|)  
< \varepsilon^2 \max(|\Delta_{n-6}|,|\Delta_{n-5}|,|\Delta_{n-4}|) 
< \cdots
\]  
\[ 
\cdots < \varepsilon^{\nu} \max( |\Delta_{n-3 \nu}|, |\Delta_{n-3 \nu+1}|, |\Delta_{n-3 \nu+2}| )
\] 
with $n-3\nu \geq N$ ($\nu$ is a natural number).  
Since $|\Delta_n|<\{ \alpha_0 \}$ and $|\Delta_n^{\prime}|<1 $ for all $n \geq 0$, 
we have 
$$
|\Delta_n| < \{ \alpha_0 \} \varepsilon^e 
\ \ {\rm and}\ \ 
|\Delta_n^{\prime}| < \varepsilon^e,\ \ {\rm where}\ \ e=(n-N)/3, 
$$
for all $n \geq N+3$. 
Letting  $c=\ve^{\frac{1}{3}}$, we have proved (I).  

(II) \ 
From Lemma 2.12,  
$\ \displaystyle\limsup_{n \to \infty} \displaystyle{\alpha_n  \over  \beta_n} <1$  
and  
$\ \displaystyle\liminf_{n \to \infty} \beta_n >1$. 
Let $M$ be a real number with $\beta_n \leq M$ for all $n \geq 0$. 
From the equality (1.2),  $r_n < 3Mr_{n-1}$ so $r_n / r_{n-1} <3M$.  
We define $a=-\log c / \log 3M \ (>0)$, i.e. $c(3M)^a=1$. 
Since 
$c^n r_n^a /c^{n-1} r_{n-1}^a =c(r_n / r_{n-1})^a < c(3M)^a=1$, 
we have 
$c^n r_n^a < c^{n-1} r_{n-1}^a < \cdots <cr_1^a =cb_1^a < cM^a =1/3^a <1$. 
Therefore, 
$c^n < 1/r_n^a$ for all $n \geq 0$ and then (II) is proved by (I).  
 \hfill  \qed

\bigskip

Although the boundedness of $\beta_n$ is one of the necessary conditions of the periodicity of 
the expansion, it turned out to be sufficient for   
$\Delta_n \to 0\ (n \to \infty)$.   
This theorem shows that $\displaystyle\lim_{n \to \infty} \Delta_n=0$  when 
$(\alpha_0,\,\beta_0)$ has the periodic expansion. 
The condition (I) is weaker than (II). 
Let $(\alpha_0,\,\beta_0)$ be a point in $\mathcal{D}$ with 
$a_{2i}=ki+1,\,a_{2i+1}=1,\,b_{2i}=(k+1)i+2,\,b_{2i+1}=5\ (i \geq 0)$, where $k$ is a fixed natural number.  
Then we can see that 
$\ \displaystyle\limsup_{n \to \infty} \displaystyle{\alpha_n  \over  \beta_n} = k/(k+1)$,   
$\ \displaystyle\liminf_{ n \to \infty} \beta_n \geq 5$.  
Thus, 
$\limsup_{n \to \infty} \{ \alpha_n \} \leq 1/5$,  $\limsup_{n \to \infty} \{ \beta_n \} =  k/(k+1)$ and  
$$ 
\max \Bigl( 
\displaystyle{ 1+\limsup_{n \to \infty} \{ \alpha_n \}  \over  2 },\ \limsup_{n \to \infty} \{ \beta_n \}  \Bigr) 
\leq \max(3/5,\ k/(k+1)). 
$$ 
Here  $\beta_n\ (n \geq 0)$ is not bounded from above.   

Perron(Satz\,V\ in[7]) proved a sufficient condition to have 
$\displaystyle\lim_{n \to \infty} \Delta_n=0$. 
In the case of dimension 2, Perron's condition is equal to that 
$\displaystyle{2+a_n  \over  b_n} \leq \theta <1$ for all sufficiently large $n$ 
where $\theta$ is a real number. 
Under this condition, it holds $2+a_n< b_n$, so $3 \leq b_n$ and 
$ \displaystyle\liminf_{n \to \infty} \beta_n \geq 3$. 
Further $\displaystyle{\alpha_n  \over  \beta_n}<\displaystyle{1+a_n  \over  b_n}< \theta$, 
thus 
$\ \displaystyle\limsup_{n \to \infty} \displaystyle{\alpha_n  \over  \beta_n}<1$. 
Therefore the condition of Theorem 2.3 includes Perron's.  
He also gave an example of $\displaystyle\lim_{n \to \infty} \Delta_n \not=0$ ($\S$6 in [7]). 
To be precise, in our notation, that example satisfies 
$\displaystyle\lim_{n \to \infty} \Delta_{2n}\not=0$, 
$\Delta_{2n} \Delta_{2n+1}<0$ and 
$\displaystyle\prod_{n \geq 0}\{ \alpha_{2n} \}\not=0$. 
Further, $\displaystyle\lim_{n \to \infty} \Delta_{2n+1}=0$, \ 
$\displaystyle\lim_{n \to \infty} \displaystyle{\alpha_{2n}  \over  \beta_{2n}}=0$ and 
$ \displaystyle\liminf_{n \to \infty} \beta_n =1$ can be seen easily. 
 
\bigskip

{\bf Lemma 2.13.}\ \ {\it If $\displaystyle\prod_{n \geq 0}\{ \beta_{n} \} \not=0$, then 
$\beta_n  \to \infty$, $\alpha_n \to \infty$ and $\{ \alpha_n \} \to 0$ as $n \to \infty$.  
}

\medskip

{\it Proof.}\ \ 
We have $\{ \beta_n \} \to 1\ ( n \to \infty)$ that is $ \alpha_n/\beta_n \to 1$. 
If $\beta_n\ (n \geq 0)$ is bounded from above, then 
$\displaystyle\prod_{n \geq 0}\{ \beta_{n} \} =0$  
from Lemma 2.12(i).  
So, $\beta_n\ (n \geq 0)$ is not bounded from above and $\alpha_n\ (n \geq 0)$ is not either.  
If $\{ \alpha_n \}$ does not converge at $0\ (n \to 0)$, then there exists a subsequence  
$\{ \alpha_{n_j} \}$ such that $\{ \alpha_{n_j} \} \to \ve \ (j \to \infty)$ for some $\ve \ (0< \ve \leq 1)$.  
We have 
$
\beta_{n_j+1}=1/\{ \alpha_{n_j} \} \to 1/\ve <\infty\ (j \to \infty)
$ 
which contradicts that 
$
\beta_n \to \infty\ (n \to \infty).
$ 
 
\hfill  \qed

\bigskip 

{\bf Theorem 2.4.} \ \ {\it                              
Assume there exist infinite number of $n$ with $\Delta_{n-1} \Delta_{n}>0$ 
and $\displaystyle\prod_{n \geq 0}\{ \beta_{n} \} \not=0$.  
Then, $\Delta_n \to 0\ (n \to \infty)$. 
}

\medskip 

{\it Proof.}\ \ 
By means of Lemma 2.13, we have  $\displaystyle\liminf_{n \to \infty} \beta_n >1$. 
So, Theorem 2.2 implies our result. 
 \hfill  \qed

\bigskip

{\bf Theorem 2.5.} \ \ {\it Assume  there are only finite number of $n$ with 
$\Delta_{n-1} \Delta_{n}$ $>0$ and  $\displaystyle\prod_{n \geq 0}\{ \beta_{n} \}=0$. 
Then, $\Delta_{n} \to 0\ (n \to \infty)$. }

\medskip

{\it Proof.}\ \ 
From the former assumption, there exists a natural number $N$ such that 
$\Delta_{n-1} \Delta_{n}<0$ for all $n \geq N$. Then, by Lemma 2.2, 
\[| \Delta_{N+m} |< \{ \beta_{N+m} \}| \Delta_{N+m-1} |<
\{ \beta_{N+m} \} \{ \beta_{N+m-1} \}| \Delta_{N+m-2} | < \cdots
\]  
\[\cdots < 
\{ \beta_{N+m} \} \{ \beta_{N+m-1} \} \{ \beta_{N+m-2} \} \cdots \{ \beta_{N+1} \} | \Delta_{N} |
\]
for all $m \geq 1$.  
By the latter assumption, the coefficient of $| \Delta_{N} |$ tends to zero as $m \to \infty$   
and the result is proved. 
\hfill  \qed

\bigskip

For integers $k$ and $\ell\,( -2 \leq k < \ell)$, we define 
\[ 
D_{k,\ell} 
= 
\left|
   \begin{array}{ccc}
     1            &   r_{k}  &  r_{\ell}   \\
     \alpha_0  &   p_{k}  &  p_{\ell}   \\
     \beta_0   &   q_{k}  &  q_{\ell}  
    \end{array} 
\right|.     
\]  
It holds that 
\[ 
D_{-2,\, -1 } 
= 
\left|
   \begin{array}{ccc}
     1            &   0  &  0   \\
     \alpha_0  &   1  &  0   \\
     \beta_0   &   0  &  1  
    \end{array} 
\right| = 1\ \ {\rm and}\ \      
D_{-1,\,0} 
= 
\left|
   \begin{array}{ccc}
     1            &   0  &  1   \\
     \alpha_0  &   0  &  a_0   \\
     \beta_0   &   1  &  b_0  
    \end{array} 
\right| = \alpha_0-a_0 = \{ \alpha_0 \}.      
\]  
By the second equality of (1.1),  letting  $\pi_{n} = \beta_1 \cdots \beta_{n}$, 
we have 
\[ 
\alpha_n 
\!\!= \!\! 
\left| \!\!         
   \begin{array}{ccc}
     r_{n-3}  &   \pi_n                &      r_{n-1}                  \\
     p_{n-3}  &   \pi_n \alpha_0    &      p_{n-1}                \\
     q_{n-3}  &   \pi_n \beta_0     &      q_{n-1}    
    \end{array} 
\!\!\right| 
\!\!=\!\! -\pi_n D_{n-3,\,n-1}  
\ {\rm and}\ 
\beta_n 
\!\!=\!\!  
\left|  \!\!        
   \begin{array}{ccc}
     r_{n-3}  &    r_{n-2}       &     \pi_n                               \\
     p_{n-3}  &    p_{n-2}       &     \pi_n \alpha_0                   \\
     q_{n-3}  &    q_{n-2}       &     \pi_n \beta_0       
    \end{array} 
\!\!\right| 
\!\!= \!\!\pi_n D_{n-3,\,n-2}.   
\]  
Thus, 
\[ 
D_{n-3,\,n-1} 
= 
-\displaystyle{ 1  \over  \beta_1 \cdots \beta_{n-1} }
\displaystyle{ \alpha_n  \over  \beta_n } 
= 
-\{ \alpha_0 \} \cdots \{ \alpha_{n-2} \} \{ \beta_{n-1} \}, 
\]  
\[ 
D_{n-3,\,n-2} 
= 
\displaystyle{ 1  \over  \beta_1 \cdots \beta_{n-1} } 
= 
\{ \alpha_0 \} \cdots \{ \alpha_{n-2} \}. 
\] 
Summing up above,  we get  the following 

\bigskip
{\bf  Lemma 2.14.}\ \ {\it For $n \geq1$, we have 
\[
D_{n-1,\, n} 
= 
\{ \alpha_0 \} \cdots \{ \alpha_{n} \}\ \ and\ \ 
D_{n-2,\, n} 
= 
-\{ \alpha_0 \} \cdots \{ \alpha_{n-1} \} \{ \beta_n \}. 
\] 
\[ 
\alpha_n 
= -\displaystyle{ D_{n-3,\, n-1}  \over  D_{n-2,\,n-1}  },\ \ 
\beta_n 
= \displaystyle{ D_{n-3,\, n-2}  \over  D_{n-2,\,n-1}  },\ \ 
\{ \alpha_n \} 
= 
\displaystyle{ D_{n-1,\, n}  \over  D_{n-2,\,n-1}  }\ \ and\ \ 
\{ \beta_n \} 
= 
-\displaystyle{ D_{n-2,\ n}  \over  D_{n-2,\,n-1}  }.\  
\]
}

\bigskip

See also  equations (6) in [7] and (1.15) in [1]. 

\bigskip

{\bf Lemma 2.15.}\ \ {\it 
It holds that 

{\rm (i)} $\Delta_n \Delta_{n+1}^{\prime}- \Delta_{n+1} \Delta_n^{\prime}>0$ \ and \ 
{\rm (ii)} $\Delta_n \Delta_{n+2}^{\prime}-\Delta_{n+2} \Delta_n^{\prime} <0$ for all $n \geq 0$. }

\medskip

{\it Proof.}\ \ 
\[ 
\Delta_n \Delta_{n+1}^{\prime}- \Delta_{n+1} \Delta_n^{\prime} 
= 
\left|
   \begin{array}{ccc}
     1    &   r_{n}                         &  r_{n+1}   \\
     0    &   \Delta_{n}                  &  \Delta_{n+1}   \\
     0    &   \Delta^{\prime}_{n}      &  \Delta^{\prime}_{n+1}  
    \end{array} 
\right| 
= 
\left|
   \begin{array}{ccc}
     1            &   r_{n}  &  r_{n+1}   \\
     \alpha_0  &   p_{n}  &  p_{n+1}   \\
     \beta_0   &   q_{n}  &  q_{n+1}  
    \end{array} 
\right| 
\] 
\[ 
= 
D_{n,\,n+1}  
= 
\{ \alpha_0 \} \cdots \{ \alpha_{n+1} \} > 0,  
\]  
\[ 
\Delta_n \Delta_{n+2}^{\prime}- \Delta_{n+2} \Delta_n^{\prime} 
= 
\left|
   \begin{array}{ccc}
     1    &   r_{n}                         &  r_{n+2}   \\
     0    &   \Delta_{n}                  &  \Delta_{n+2}   \\
     0    &   \Delta^{\prime}_{n}      &  \Delta^{\prime}_{n+2}  
    \end{array} 
\right| 
= 
\left|
   \begin{array}{ccc}
     1            &   r_{n}  &  r_{n+2}   \\
     \alpha_0  &   p_{n}  &  p_{n+2}   \\
     \beta_0   &   q_{n}  &  q_{n+2}  
    \end{array} 
\right| 
\] 
\[ 
= 
D_{n,\,n+2}  
= 
-\{ \alpha_0 \} \cdots \{ \alpha_{n+1} \} \{ \beta_{n+2} \} < 0.   
\] 
\hfill  \qed

\bigskip

{\bf Lemma 2.16.} \ \ {\it The number of $n$ with $\Delta_{n-1} \Delta_n >0$ is finite 
if and only if the number of $n$ with $\Delta_{n-1}^{\prime} \Delta_n^{\prime} >0$ is finite. }

\medskip

{\it Proof.}\ \ Sufficiency ; let $N$ be a number such that 
$\Delta_{n-1} \Delta_{n}<0$ for all $n \geq N$. 
Suppose $\Delta_{m}>0$ and $\Delta_{m}^{\prime}<0$ for some $m \geq N$, 
then from Lemma 2.15(i), we see $\Delta_{m+1}^{\prime}>0$ 
because $\Delta_{m+1}<0$,   
and from  Lemma 2.15(ii), $\Delta_{m+2}^{\prime}<0$ because $\Delta_{m+2}>0$.  
Repeating this, we have 
$\Delta_{n}^{\prime} \Delta_{n+1}^{\prime}<0$ for all $n \geq m$. 
If there is no number $m (\geq N)$ with   
$\Delta_{m}>0$ and $\Delta_{m}^{\prime}<0$, then 
$\Delta_{m}^{\prime}>0$, $\Delta_{m+2}^{\prime}>0$, $\Delta_{m+4}^{\prime}>0 \cdots$ 
(i.e., $\Delta_m^{\prime}>0$ whenever $\Delta_m >0$). 
Therefore $\Delta_{n}^{\prime} \Delta_{n+1}^{\prime}<0$ for all $n \geq m$. 

Necessity ; let $N^{\prime}$ be a number such that 
$\Delta_{n-1}^{\prime} \Delta_{n}^{\prime}<0$ for all $n \geq N^{\prime}$. 
Suppose $\Delta_{m}^{\prime}>0$ and $\Delta_{m}>0$ for some $m \geq N^{\prime}$, 
then from  Lemma 2.15(i), we see $\Delta_{m+1}<0$ 
because $\Delta_{m+1}^{\prime}<0$,   
and from  Lemma 2.15(ii) $\Delta_{m+2}>0$ because $\Delta_{m+2}^{\prime}>0$.  
Repeating this, we have 
$\Delta_{n} \Delta_{n+1}<0$ for all $n \geq m$. 
If there is no number $m (\geq N^{\prime})$ with   
$\Delta_{m}^{\prime}>0$ and $\Delta_{m}>0$, then 
$\Delta_{m}<0$, $\Delta_{m+2}<0$, $\Delta_{m+4}<0 \cdots$ 
(i.e., $\Delta_m <0$ whenever $\Delta_m^{\prime} >0$). 
Thus $\Delta_{n} \Delta_{n+1}<0$ for all $n \geq m$ and the number of $n$ with  
$\Delta_{n-1} \Delta_n >0$ is finite. 
\hfill  \qed

\bigskip

Last, we shall give a remark on the case where $\beta_n \to \infty\ (n \to \infty)$.   
Because  
$ \beta_n \Delta_{n-1}+\alpha_n \Delta_{n-2}+\Delta_{n-3}=0$ from  (1.3),  
we see  $\Delta_n \to 0\ (n \to \infty)$ when 
$\ \displaystyle{  \alpha_{n}+1  \over  \beta_{n}  }  \to 0\ (n \to \infty)$.  
In  the case of Perron's example  with $\lim_{n \to \infty} \Delta_n \not= 0$ $(\S6\ {\rm in}\ [7])$ mentioned before,  
$\ \displaystyle{  \alpha_{n}  \over  \beta_{n}  }  \to 0\ (n \to \infty)$  
but 
$ \displaystyle{  1  \over  \beta_{n}  }$ does not converge at 0. 

\vspace{1cm}

{\bf 3.}\  Conjugates of $\alpha_n$ and $\beta_n$.      

\medskip
 
Let $\gamma$ be a real algebraic irrational. 
The classical (i.e., of dimension 1) continued fraction of $\gamma=\gamma_0$ is defined by
$\gamma_{n-1}=[\gamma_{n-1}]+\displaystyle{1  \over  \gamma_n}\ \ ( \gamma_n >1,\ n \geq 1).$ 
Let $p_n/r_n\ (n \geq 0)$ denote the convergent of $\gamma$, 
then 
$ \gamma=\displaystyle{ p_n \gamma_n+p_{n-1}  \over  r_n \gamma_n+r_{n-1}}$.  
We denote by $\gamma^{\prime}$ an algebraic conjugate of  $\gamma$.  
Let   
$ \gamma_n^{\prime}=-\displaystyle{  \gamma^{\prime}r_{n-1}-p_{n-1}  \over  \gamma^{\prime} r_n-p_{n}}$,  
then  
$\gamma_n^{\prime}\ $ is an  algebraic conjugate of $\gamma_n$. 
From this,  
$$
\gamma_n^{\prime}+\displaystyle{r_{n-1}  \over  r_n} =\displaystyle{r_{n-1}  \over  r_n} 
\displaystyle{ p_{n} / r_{n} - p_{n-1} / r_{n-1}   \over  p_{n} / r_{n} - \gamma^{\prime} }
$$ 
so                                                             

$\hspace{45mm} \displaystyle{\lim_{n \to \infty}}  \Bigl( \gamma_n^{\prime}+\displaystyle{r_{n-1}  \over  r_n}  \Bigr)=0.   \hfill    (3.1)  $

\medskip
\noindent
Further, when $\gamma^{\prime}$ is real we can see $-1 < \gamma_n^{\prime} < 0$ for all sufficiently large $n$. 
How is in the case of the Jacobi-Perron algorithm of dimension 2? 

Let $\alpha_0$ and $\beta_0$ be real algebraic irrationals. 
We denote by $\alpha_0^{\prime}$ $(\not=  \alpha_0)$ and $\beta_0^{\prime}$ $(\not=  \beta_0)$ 
an algebraic conjugate of  
$\alpha_0$ and $\beta_0$, respectively. 
Let $\alpha_n^{\prime}$ and $\beta_n^{\prime}\ (n \geq 1)$  
be algebraic numbers defined by the relations 
 $\alpha_{n-1}^{\prime}=a_{n-1}+\displaystyle{ 1 \over \beta_n^{\prime}}$ and 
$\beta_{n-1}^{\prime}=b_{n-1}+\displaystyle{ \alpha_n^{\prime} \over \beta_n^{\prime}}\ \ (n \geq 1)$  
successively. 
Then $\alpha_n^{\prime}$ and $\beta_n^{\prime}$ are algebraic conjugates of 
$\alpha_n$ and $\beta_n$.   

Define  
\[ 
D_{k,\,\ell}^{\prime}
= 
\left|
   \begin{array}{ccc}
     1            &   r_{k}  &  r_{\ell}   \\
     \alpha_0^{\prime}  &   p_{k}  &  p_{\ell}   \\
     \beta^{\prime}_0   &   q_{k}  &  q_{\ell}  
    \end{array} 
\right|.
\]  
Following the classical case (3.1), the behavior of 
$
\gamma_n^{\prime}+\displaystyle{r_{n-1}  \over  r_n} 
$ 
according to $n \to \infty$,  
we shall study about 
$
\beta_n^{\prime}+\displaystyle{r_{n-2}  \over r_{n-1} }\alpha_n^{\prime}
+\displaystyle{r_{n-3}  \over r_{n-1} }
$. 
From Lemma 2.14, it holds 

$
\beta_n^{\prime}+\displaystyle{r_{n-2}  \over r_{n-1} }\alpha_n^{\prime}
+\displaystyle{r_{n-3}  \over r_{n-1} } 
= 
\displaystyle{ D_{n-3,\,n-2}^{\prime}  \over  D_{n-2,\,n-1}^{\prime} } 
- 
\displaystyle{r_{n-2}  \over r_{n-1} } 
\displaystyle{ D_{n-3,\,n-1}^{\prime}  \over  D_{n-2,\,n-1}^{\prime} } 
+ 
\displaystyle{r_{n-3}  \over r_{n-1} } 
$

\medskip
 
$ \hspace{30mm}
=
\displaystyle{ r_{n-1} D_{n-3,\,n-2}^{\prime} -  r_{n-2} D_{n-3,\,n-1}^{\prime} + r_{n-3} D_{n-2,\,n-1}^{\prime}    
\over  r_{n-1} D_{n-2,\,n-1}^{\prime} }.     \hfill(3.2)
$                                                    

\medskip
\noindent
Let  $N_u$  and $D_e$ denote, respectively, the numerator  and the denominator of $(3.2)$.  
\[ 
N_u=r_{n-1} 
\left|
   \begin{array}{ccc}
     1                        &   r_{n-3}  &  r_{n-2}   \\
     \alpha_0^{\prime}  &   p_{n-3}  &  p_{n-2}   \\
     \beta^{\prime}_0   &   q_{n-3}  &  q_{n-2}  
    \end{array} 
\right| 
-r_{n-2} 
\left|
   \begin{array}{ccc}
     1                        &   r_{n-3}  &  r_{n-1}   \\
     \alpha_0^{\prime}  &   p_{n-3}  &  p_{n-1}   \\
     \beta^{\prime}_0   &   q_{n-3}  &  q_{n-1}  
    \end{array} 
\right| 
+r_{n-3} 
\left|
   \begin{array}{ccc}
     1                        &   r_{n-2}  &  r_{n-1}   \\
     \alpha_0^{\prime}  &   p_{n-2}  &  p_{n-1}   \\
     \beta^{\prime}_0   &   q_{n-2}  &  q_{n-1}  
    \end{array} 
\right|
\] 

\noindent
$= 
\left|
   \begin{array}{ccc}
     1                        &   r_{n-3}  &   0   \\
     \alpha_0^{\prime}  &   p_{n-3}  &   r_{n-1}p_{n-2}- r_{n-2}p_{n-1}  \\
     \beta^{\prime}_0   &   q_{n-3}  &   r_{n-1}q_{n-2}- r_{n-2}q_{n-1} 
    \end{array} 
\right| 
+ 
\displaystyle{r_{n-3}  \over r_{n-1} } 
\left|
   \begin{array}{ccc}
     1                        &  r_{n-1} r_{n-2}  &  r_{n-1}   \\
     \alpha_0^{\prime}  &   r_{n-1}p_{n-2}  &  p_{n-1}   \\
     \beta^{\prime}_0   &  r_{n-1} q_{n-2}  &  q_{n-1}  
    \end{array} 
\right| 
$

\medskip
\noindent
$= 
-r_{n-3}
\left|
   \begin{array}{ccc}
     1                        &   0                                      &   1   \\
     \alpha_0^{\prime}  &   \!r_{n-1}p_{n-2}- r_{n-2}p_{n-1}  &  \!\displaystyle{ p_{n-3}  \over  r_{n-3} } \\
     \beta^{\prime}_0   &  \! r_{n-1}q_{n-2}- r_{n-2}q_{n-1}  &  \!\displaystyle{ q_{n-3}  \over  r_{n-3} }
    \end{array} 
\right| 
+ 
\displaystyle{r_{n-3}  \over r_{n-1} } 
\left|
   \begin{array}{ccc}
     1                        &    0                                        &   \!r_{n-1}   \\
     \alpha_0^{\prime}  &    \!r_{n-1}p_{n-2}- r_{n-2}p_{n-1}    &   \!p_{n-1}   \\
     \beta^{\prime}_0   &    \!r_{n-1}q_{n-2}- r_{n-2}q_{n-1}    &  \!q_{n-1}  
    \end{array} 
\right| 
$ 

\bigskip
\noindent
$= 
-r_{n-3}
\left|
   \begin{array}{ccc}
     1                        &   0                                      &   1   \\
     \alpha_0^{\prime}  &   \!r_{n-1}p_{n-2}- r_{n-2}p_{n-1}  &  \!\displaystyle{ p_{n-3}  \over  r_{n-3} } \\
     \beta^{\prime}_0   &  \! r_{n-1}q_{n-2}- r_{n-2}q_{n-1}  &  \!\displaystyle{ q_{n-3}  \over  r_{n-3} }
    \end{array} 
\right| 
+ 
r_{n-3}  
\left|
   \begin{array}{ccc}
     1                        &    0                                        &   1   \\
     \alpha_0^{\prime}  &    \!r_{n-1}p_{n-2}- r_{n-2}p_{n-1}    &   \!\displaystyle{ p_{n-1}  \over  r_{n-1} }   \\
     \beta^{\prime}_0   &    \!r_{n-1}q_{n-2}- r_{n-2}q_{n-1}    &  \!\displaystyle{ q_{n-1}  \over  r_{n-1} }   
    \end{array} 
\right| 
$ 

\bigskip
\noindent
$= 
r_{n-3}  
\left|
   \begin{array}{ccc}
 1                        &    0                                        &   0   \\
 \alpha_0^{\prime}  &    \!r_{n-1}p_{n-2}- r_{n-2}p_{n-1}    &   \!\displaystyle{ p_{n-1}  \over  r_{n-1} } -  \displaystyle{ p_{n-3}  \over  r_{n-3} } \\
 \beta^{\prime}_0   &    \!r_{n-1}q_{n-2}- r_{n-2}q_{n-1}    &  \!\displaystyle{ q_{n-1}  \over  r_{n-1} } -  \displaystyle{ q_{n-3}  \over  r_{n-3} }  
    \end{array} 
\right| 
$

\bigskip
\noindent
$= 
r_{n-3}  
\left|
   \begin{array}{ccc}
 1                        &    0                                        &   0   \\
 \alpha_0^{\prime}  &    \!r_{n-1}\Delta_{n-2}- r_{n-2}\Delta_{n-1}    &   \!\displaystyle{ p_{n-1}  \over  r_{n-1} } -  \displaystyle{ p_{n-3}  \over  r_{n-3} } \\
 \beta^{\prime}_0   &   \!r_{n-1}\Delta_{n-2}^{\prime}- r_{n-2}\Delta_{n-1}^{\prime}  &  \!\displaystyle{ q_{n-1}  \over  r_{n-1} } -  \displaystyle{ q_{n-3}  \over  r_{n-3} }  
    \end{array} 
\right| 
$

\bigskip
\noindent
$= 
r_{n-3}  
\Big\{ 
(r_{n-1}\Delta_{n-2}- r_{n-2}\Delta_{n-1}) \Big(\displaystyle{ q_{n-1}  \over  r_{n-1} } \!-  \!\displaystyle{ q_{n-3}  \over  r_{n-3} } \Big) 
\!-\! 
(r_{n-1}\Delta_{n-2}^{\prime}- r_{n-2}\Delta_{n-1}^{\prime}) \Big(\displaystyle{ p_{n-1}  \over  r_{n-1} } \!- \! \displaystyle{ p_{n-3}  \over  r_{n-3} } \Big) 
\Big\}.
$

\bigskip
\noindent
Therefore, 

\medskip
\noindent
$
|N_u|< 
r_{n-3}  
\Big\{ 
(r_{n-1}|\Delta_{n-2}|+ r_{n-2}|\Delta_{n-1}|) \Big|\displaystyle{ q_{n-1}  \over  r_{n-1} } \!-\!  \displaystyle{ q_{n-3}  \over  r_{n-3} } \Big| 
$    

\medskip
  
$
\hspace{30mm}+ 
(r_{n-1}|\Delta_{n-2}^{\prime}|+ r_{n-2}|\Delta_{n-1}^{\prime}|) \Big|\displaystyle{ p_{n-1}  \over  r_{n-1} } \!-\!  \displaystyle{ p_{n-3}  \over  r_{n-3} } \Big| 
\Big\}.        
$      \hfill(3.3)                                        

\medskip
$ 
D_e 
= 
r_{n-1} 
\left|
   \begin{array}{ccc}
     1                        &   r_{n-2}  &  r_{n-1}   \\
     \alpha_0^{\prime}  &   p_{n-2}  &  p_{n-1}   \\
     \beta^{\prime}_0   &   q_{n-2}  &  q_{n-1}  
    \end{array} 
\right| 
= 
\left|
   \begin{array}{ccc}
     1                        &   r_{n-1} r_{n-2}  &  r_{n-1}   \\
     \alpha_0^{\prime}  &   r_{n-1} p_{n-2}  &  p_{n-1}   \\
     \beta^{\prime}_0   &   r_{n-1} q_{n-2}  &  q_{n-1}  
    \end{array} 
\right| 
$ 

\medskip
\noindent
$ 
= 
\left|
   \begin{array}{ccc}
     1                        &    0                                        &   \!r_{n-1}   \\
     \alpha_0^{\prime}  &    \!r_{n-1}p_{n-2}- r_{n-2}p_{n-1}    &   \!p_{n-1}   \\
     \beta^{\prime}_0   &    \!r_{n-1}q_{n-2}- r_{n-2}q_{n-1}    &  \!q_{n-1}  
    \end{array} 
\right| 
= 
r_{n-1}
\left|
   \begin{array}{ccc}
     1                        &    0                                          &   1   \\
     \alpha_0^{\prime}  &    \!r_{n-1}p_{n-2}- r_{n-2}p_{n-1}    &   \!\displaystyle{ p_{n-1}  \over  r_{n-1} }   \\
     \beta^{\prime}_0   &    \!r_{n-1}q_{n-2}- r_{n-2}q_{n-1}    &  \!\displaystyle{ q_{n-1}  \over  r_{n-1} }   
    \end{array} 
\right| 
$ 

\medskip
\noindent
$ 
= 
r_{n-1}
\left|
   \begin{array}{ccc}
     1                        &    0                                          &   0                          \\
     \alpha_0^{\prime}  &    \!r_{n-1}p_{n-2}- r_{n-2}p_{n-1}    &   \!\displaystyle{ p_{n-1}  \over  r_{n-1} } - \alpha_0^{\prime}     \\
     \beta^{\prime}_0   &    \!r_{n-1}q_{n-2}- r_{n-2}q_{n-1}    &  \!\displaystyle{ q_{n-1}  \over  r_{n-1} } - \beta^{\prime}_0  
    \end{array} 
\right| 
$

\medskip
\noindent
$ 
=
r_{n-1}  
\Big\{ 
(r_{n-1}\Delta_{n-2}- r_{n-2}\Delta_{n-1}) \Big(\displaystyle{ q_{n-1}  \over  r_{n-1} } \!-  \beta^{\prime}_0 \Big) 
$ 

\noindent
$
\hspace{40mm}
- 
(r_{n-1}\Delta_{n-2}^{\prime}- r_{n-2}\Delta_{n-1}^{\prime}) \Big(\displaystyle{ p_{n-1}  \over  r_{n-1} } \!- \! \alpha_0^{\prime}  \Big) 
\Big\}.
$ \hfill(3.4)                                  

\medskip

Suppose 
the number of  \,$n\,( \geq 1)$ with $\Delta_{n-1} \Delta_n >0$ is finite. 
Let $N$ be a natural number such that   $\Delta_{n-1} \Delta_n <0$ and $\Delta_{n-1}^{\prime} \Delta_n^{\prime} <0$ 
for all $n>N$ (see Lemma 2.16).  
We shall consider the following four cases;

(1)\   $\Delta_N \Delta_N^{\prime}<0$ and $\Delta_{n-2} >0\ (n-2 \geq N)$ for some $n$. 
Because  
$\Delta_{n-1} <0$, $\Delta_{n-2}^{\prime} <0$ and $\Delta_{n-1}^{\prime} >0$, from (3.4), 

\medskip

$
D_e=   
r_{n-1}  
\Big\{ 
(r_{n-1}|\Delta_{n-2}|+ r_{n-2}|\Delta_{n-1}|) \Big( \displaystyle{ q_{n-1}  \over  r_{n-1} } \!-  \beta^{\prime}_0 \Big) 
$ 
 
$
\hspace{60mm}
+ 
(r_{n-1}|\Delta_{n-2}^{\prime}|+ r_{n-2}|\Delta_{n-1}^{\prime}|) \Big(\displaystyle{ p_{n-1}  \over  r_{n-1} } \!- \! \alpha_0^{\prime}  \Big) 
\Big\}.
$   

(2)\  $\Delta_N \Delta_N^{\prime}<0$ and $\Delta_{n-2} <0\ (n-2 \geq N)$ for some $n$. 
Because  
$\Delta_{n-1} >0$, $\Delta_{n-2}^{\prime} >0$ and $\Delta_{n-1}^{\prime} <0$, from (3.4), 

\medskip

$
D_e=   
-r_{n-1}  
\Big\{ 
(r_{n-1}|\Delta_{n-2}|+ r_{n-2}|\Delta_{n-1}|) \Big( \displaystyle{ q_{n-1}  \over  r_{n-1} } \!-  \beta^{\prime}_0 \Big) 
$ 
 
$
\hspace{60mm}
+ 
(r_{n-1}|\Delta_{n-2}^{\prime}|+ r_{n-2}|\Delta_{n-1}^{\prime}|) \Big(\displaystyle{ p_{n-1}  \over  r_{n-1} } \!- \! \alpha_0^{\prime}  \Big) 
\Big\}.
$   

(3)\  $\Delta_N \Delta_N^{\prime}>0$ and $\Delta_{n-2} >0\ (n-2 \geq N)$ for some $n$. 
Because  
$\Delta_{n-1} <0$, $\Delta_{n-2}^{\prime} >0$ and $\Delta_{n-1}^{\prime} <0$, from (3.4), 

\medskip

$
D_e=   
r_{n-1}  
\Big\{ 
(r_{n-1}|\Delta_{n-2}|+ r_{n-2}|\Delta_{n-1}|) \Big( \displaystyle{ q_{n-1}  \over  r_{n-1} } \!-  \beta^{\prime}_0 \Big) 
$ 
 
$
\hspace{60mm}
- 
(r_{n-1}|\Delta_{n-2}^{\prime}|+ r_{n-2}|\Delta_{n-1}^{\prime}|) \Big(\displaystyle{ p_{n-1}  \over  r_{n-1} } \!- \! \alpha_0^{\prime}  \Big) 
\Big\}.
$   

(4)\  $\Delta_N \Delta_N^{\prime}>0$ and $\Delta_{n-2} <0\ (n-2 \geq N)$ for some $n$. 
Because  
$\Delta_{n-1} >0$, $\Delta_{n-2}^{\prime} <0$ and $\Delta_{n-1}^{\prime} >0$, from (3.4), 

\medskip

$
D_e=   
-r_{n-1}  
\Big\{ 
(r_{n-1}|\Delta_{n-2}|+ r_{n-2}|\Delta_{n-1}|) \Big( \displaystyle{ q_{n-1}  \over  r_{n-1} } \!-  \beta^{\prime}_0 \Big) 
$ 
 
$
\hspace{60mm}
- 
(r_{n-1}|\Delta_{n-2}^{\prime}|+ r_{n-2}|\Delta_{n-1}^{\prime}|) \Big(\displaystyle{ p_{n-1}  \over  r_{n-1} } \!- \! \alpha_0^{\prime}  \Big) 
\Big\}.
$   

\noindent
Therefore, when $\Delta_N \Delta_N^{\prime}<0$ it holds   

\medskip

$
|D_e|=   
r_{n-1}  
\Big| 
(r_{n-1}|\Delta_{n-2}|+ r_{n-2}|\Delta_{n-1}|) \Big( \displaystyle{ q_{n-1}  \over  r_{n-1} } \!-  \beta^{\prime}_0 \Big) 
$ 
 
$
\hspace{55mm}
+ 
(r_{n-1}|\Delta_{n-2}^{\prime}|+ r_{n-2}|\Delta_{n-1}^{\prime}|) \Big(\displaystyle{ p_{n-1}  \over  r_{n-1} } \!- \! \alpha_0^{\prime}  \Big) 
\Big|,
$   \hfill(3.5)                                   

\noindent
and when $\Delta_N \Delta_N^{\prime}>0$,  

\medskip

$
|D_e|=   
r_{n-1}  
\Big| 
(r_{n-1}|\Delta_{n-2}|+ r_{n-2}|\Delta_{n-1}|) \Big( \displaystyle{ q_{n-1}  \over  r_{n-1} } \!-  \beta^{\prime}_0 \Big) 
$ 
 
$
\hspace{55mm}
- 
(r_{n-1}|\Delta_{n-2}^{\prime}|+ r_{n-2}|\Delta_{n-1}^{\prime}|) \Big(\displaystyle{ p_{n-1}  \over  r_{n-1} } \!- \! \alpha_0^{\prime}  \Big) 
\Big|. 
$   \hfill(3.6)                       

\bigskip
 
{\bf Theorem 3.1.}\ \  {\it         
Suppose the number of  \,$n\,( \geq 1)$ with $\Delta_{n-1} \Delta_n >0$ is finite  
and both $ \alpha_0^{\prime} $ and $ \beta_0^{\prime} $ are real.    
Let N be the same as before.   
If  $\Delta_N \Delta_N^{\prime} (\alpha_0 - \alpha_0^{\prime}) ( \beta_0 - \beta_0^{\prime} ) <0$, 
then }
\[ 
\lim_{n \to \infty}
\Big(\beta_n^{\prime}+\displaystyle{r_{n-2}  \over r_{n-1} }\alpha_n^{\prime}
+\displaystyle{r_{n-3}  \over r_{n-1} } \Big)=0.
\]  

{\it Proof.}\ \ If  $\Delta_N \Delta_N^{\prime}<0$, then $(\alpha_0 - \alpha_0^{\prime}) ( \beta_0 - \beta_0^{\prime} ) >0$  
and  
\[
\Big( \displaystyle{ q_{n-1}  \over  r_{n-1} } -  \beta^{\prime}_0 \Big) 
\Big(\displaystyle{ p_{n-1}  \over  r_{n-1} } -  \alpha_0^{\prime}  \Big) >0
\] 
for all sufficiently large $n$. From (3.5), 

$
|D_e|=   
r_{n-1}  
\Big\{ 
(r_{n-1}|\Delta_{n-2}|+ r_{n-2}|\Delta_{n-1}|) \Big| \displaystyle{ q_{n-1}  \over  r_{n-1} } \!-  \beta^{\prime}_0 \Big| 
$ 

$
\hspace{40mm}
+ 
(r_{n-1}|\Delta_{n-2}^{\prime}|+ r_{n-2}|\Delta_{n-1}^{\prime}|) \Big|\displaystyle{ p_{n-1}  \over  r_{n-1} } \!- \! \alpha_0^{\prime}  \Big| 
\Big\}.
$   \hfill(3.7)                       

If  $\Delta_N \Delta_N^{\prime}>0$, then $(\alpha_0 - \alpha_0^{\prime}) ( \beta_0 - \beta_0^{\prime} ) <0$ 
and  
\[
\Big( \displaystyle{ q_{n-1}  \over  r_{n-1} } -  \beta^{\prime}_0 \Big) 
\Big(\displaystyle{ p_{n-1}  \over  r_{n-1} } -  \alpha_0^{\prime}  \Big) <0
\] 
for all sufficiently large $n$. Then, from (3.6), we also have (3.7). 
To sum up, by (3.3) and (3.7),   
\[ 
\Big| \beta_n^{\prime}+\displaystyle{r_{n-2}  \over r_{n-1} }\alpha_n^{\prime}
+\displaystyle{r_{n-3}  \over r_{n-1} } \Big| 
= \Big| \displaystyle{ N_u  \over  D_e } \Big| 
\]

\noindent
$ 
< \!\!\displaystyle{ r_{n-3}  \over  r_{n-1} } 
\displaystyle{ 
(r_{n-1}|\Delta_{n-2}| \!+\! r_{n-2}|\Delta_{n-1}|) \!\Big|\displaystyle{ q_{n-1}  \over  r_{n-1} } \!-\!  \displaystyle{ q_{n-3}  \over  r_{n-3} } \!\Big| 
\!\!+\!\! 
(r_{n-1}|\Delta_{n-2}^{\prime}| \!+\! r_{n-2}|\Delta_{n-1}^{\prime}|) \!\Big|\displaystyle{ p_{n-1}  \over  r_{n-1} } \!-\!  \displaystyle{ p_{n-3}  \over  r_{n-3} } \!\Big|  
\over   
(r_{n-1}|\Delta_{n-2}| \!+\! r_{n-2}|\Delta_{n-1}|) \Big| \displaystyle{ q_{n-1}  \over  r_{n-1} } \!-\!  \beta^{\prime}_0 \Big| 
\!+\! 
(r_{n-1}|\Delta_{n-2}^{\prime}| \!+\! r_{n-2}|\Delta_{n-1}^{\prime}|) \Big|\displaystyle{ p_{n-1}  \over  r_{n-1} } \!- \! \alpha_0^{\prime}  \Big|   
}.  
$

For any $\varepsilon>0$,  inequalities
$
\displaystyle{ 
\Big| \displaystyle{ q_{n-1}  \over  r_{n-1} } -  \displaystyle{ q_{n-3}  \over  r_{n-3} } \Big|  
\over  
\Big| \displaystyle{ q_{n-1}  \over  r_{n-1} } - \beta^{\prime}_0 \Big| 
} 
< \varepsilon
$ 
and 
$
\displaystyle{ 
\Big| \displaystyle{ p_{n-1}  \over  r_{n-1} } -  \displaystyle{ p_{n-3}  \over  r_{n-3} } \Big|  
\over  
\Big| \displaystyle{ p_{n-1}  \over  r_{n-1} } - \alpha^{\prime}_0 \Big| 
} 
< \varepsilon
$ 
hold for all sufficiently large $n$.  Then,  

\medskip
$
\Big| \beta_n^{\prime}+\displaystyle{r_{n-2}  \over r_{n-1} }\alpha_n^{\prime}
+\displaystyle{r_{n-3}  \over r_{n-1} } \Big| 
< 
\displaystyle{ r_{n-3}  \over  r_{n-1} } \varepsilon 
< 
\varepsilon 
$   
 and our theorem is proved. \hfill \qed

\bigskip

When there exists another conjugate $\alpha_0^{\prime \prime}$ of $\alpha_0$ such that 
$\alpha_0^{\prime} < \alpha_0 < \alpha_0^{\prime \prime}$  or    
$\alpha_0^{\prime \prime} < \alpha_0 < \alpha_0^{\prime}$,  
we can change the sign of  
$\Delta_N \Delta_N^{\prime}(\alpha_0 - \alpha_0^{\prime}) ( \beta_0 - \beta_0^{\prime} ) $   
by changing $\alpha_0^{\prime}$ for $\alpha_0^{\prime \prime}$, if necessary,   
and just the same with  $\beta_0$.

\bigskip
 
{\bf Theorem 3.2.}\ \  {\it         
Suppose the number of  \,$n\,( \geq 1)$ with $\Delta_{n-1} \Delta_n >0$ is finite  
and both $ \alpha_0^{\prime} $ and $ \beta_0^{\prime} $ are imaginary.    
Let N be the same as before.   
If  $\Delta_N \Delta_N^{\prime} \Im(\alpha_0^{\prime})  \Im(\beta_0^{\prime} ) <0$, 
where $\Im(*)$ denotes the imaginary part of a complex number $*$, then 
}
\[ 
\lim_{n \to \infty}
\Big(\beta_n^{\prime}+\displaystyle{r_{n-2}  \over r_{n-1} }\alpha_n^{\prime}
+\displaystyle{r_{n-3}  \over r_{n-1} } \Big)=0.
\]  

{\it Proof.}\ \ If $\Im(\alpha_0^{\prime})  \Im(\beta_0^{\prime} ) > 0$, then  $\Delta_N \Delta_N^{\prime} < 0 $. 
It holds 

\noindent
$\Im \Big(
(r_{n-1}|\Delta_{n-2}|+ r_{n-2}|\Delta_{n-1}|) \Big( \displaystyle{ q_{n-1}  \over  r_{n-1} } \!-  \beta^{\prime}_0 \Big) 
+ 
(r_{n-1}|\Delta_{n-2}^{\prime}|+ r_{n-2}|\Delta_{n-1}^{\prime}|) \Big(\displaystyle{ p_{n-1}  \over  r_{n-1} } \!- \! \alpha_0^{\prime}  \Big) 
\Big)
$

\medskip
$=
(r_{n-1}|\Delta_{n-2}|+ r_{n-2}|\Delta_{n-1}|) ( -\Im (\beta^{\prime}_0 ) )
+ 
(r_{n-1}|\Delta_{n-2}^{\prime}|+ r_{n-2}|\Delta_{n-1}^{\prime}|) (-\Im (\alpha_0^{\prime} ) ).
$

\medskip
\noindent
So, 

\medskip
\noindent
$\Big| \Im \Big(
(r_{n-1}|\Delta_{n-2}|+ r_{n-2}|\Delta_{n-1}|) \Big( \displaystyle{ q_{n-1}  \over  r_{n-1} } \!-  \beta^{\prime}_0 \Big) 
+ 
(r_{n-1}|\Delta_{n-2}^{\prime}|+ r_{n-2}|\Delta_{n-1}^{\prime}|) \Big(\displaystyle{ p_{n-1}  \over  r_{n-1} } \!- \! \alpha_0^{\prime}  \Big) 
\Big) \Big|
$

\medskip
$=
(r_{n-1}|\Delta_{n-2}|+ r_{n-2}|\Delta_{n-1}|) | \Im (\beta^{\prime}_0 ) |
+ 
(r_{n-1}|\Delta_{n-2}^{\prime}|+ r_{n-2}|\Delta_{n-1}^{\prime}|) |\Im (\alpha_0^{\prime} ) |.
$ 

\medskip
\noindent
If  $\Im (\alpha_0^{\prime} ) \Im (\beta^{\prime}_0 ) <0$, then $\Delta_N \Delta_N^{\prime}>0$ and we  have, in the same way, 

\medskip
\noindent
$\Big| \Im \Big(
(r_{n-1}|\Delta_{n-2}|+ r_{n-2}|\Delta_{n-1}|) \Big( \displaystyle{ q_{n-1}  \over  r_{n-1} } \!-  \beta^{\prime}_0 \Big) 
- 
(r_{n-1}|\Delta_{n-2}^{\prime}|+ r_{n-2}|\Delta_{n-1}^{\prime}|) \Big(\displaystyle{ p_{n-1}  \over  r_{n-1} } \!- \! \alpha_0^{\prime}  \Big) 
\Big) \Big|
$

\medskip
$=
(r_{n-1}|\Delta_{n-2}|+ r_{n-2}|\Delta_{n-1}|) | \Im (\beta^{\prime}_0 ) |
+ 
(r_{n-1}|\Delta_{n-2}^{\prime}|+ r_{n-2}|\Delta_{n-1}^{\prime}|) |\Im (\alpha_0^{\prime} ) |.
$

\medskip
\noindent
In the case of  $\Delta_N \Delta_N^{\prime} <0$ and $\Im (\alpha_0^{\prime} ) \Im (\beta^{\prime}_0 ) >0$,  
from (3.5), we have  

 \medskip
 $
|D_e| \geq   
r_{n-1}  
\Big\{ 
(r_{n-1}|\Delta_{n-2}|+ r_{n-2}|\Delta_{n-1}|) | \Im (\beta^{\prime}_0 ) | 
$ 

$
\hspace{40mm}
+ 
(r_{n-1}|\Delta_{n-2}^{\prime}|+ r_{n-2}|\Delta_{n-1}^{\prime}|) | \Im (\alpha^{\prime}_0 ) | 
\Big\}  
$     \hfill(3.8)                       

 \noindent
because $|z| \geq |\Im(z)|$ for any complex number $z$.   

\noindent 
In the case of $\Delta_N \Delta_N^{\prime} >0$ and $\Im (\alpha_0^{\prime} ) \Im (\beta^{\prime}_0 ) <0$,  
from (3.6), we also have  (3.8). 
Therefore, under the assumption, (3.3) and (3.8)  shows 
\[ 
\Big| \beta_n^{\prime}+\displaystyle{r_{n-2}  \over r_{n-1} }\alpha_n^{\prime}
+\displaystyle{r_{n-3}  \over r_{n-1} } \Big| 
= \Big| \displaystyle{ N_u  \over  D_e } \Big| 
\]

\noindent
$ 
< \!\!\displaystyle{ r_{n-3}  \over  r_{n-1} } 
\displaystyle{ 
(r_{n-1}|\Delta_{n-2}| \!+\! r_{n-2}|\Delta_{n-1}|) \!\Big|\displaystyle{ q_{n-1}  \over  r_{n-1} } \!-\!  \displaystyle{ q_{n-3}  \over  r_{n-3} } \!\Big| 
\!\!+\!\! 
(r_{n-1}|\Delta_{n-2}^{\prime}| \!+\! r_{n-2}|\Delta_{n-1}^{\prime}|) \!\Big|\displaystyle{ p_{n-1}  \over  r_{n-1} } \!-\!  \displaystyle{ p_{n-3}  \over  r_{n-3} } \!\Big|  
\over   
(r_{n-1}|\Delta_{n-2}| \!+\! r_{n-2}|\Delta_{n-1}|) | \Im( \beta^{\prime}_0 ) | 
\!+\! 
(r_{n-1}|\Delta_{n-2}^{\prime}| \!+\! r_{n-2}|\Delta_{n-1}^{\prime}|) | \Im(\alpha_0^{\prime} ) |   
}.  
$

\medskip
For any $\varepsilon>0$,  
\[
\displaystyle{ 
\Big| \displaystyle{ q_{n-1}  \over  r_{n-1} } -  \displaystyle{ q_{n-3}  \over  r_{n-3} } \Big|  
\over  
| \Im (\beta^{\prime}_0) | 
} 
< \varepsilon\ \ {\rm and}\ \ 
  \displaystyle{ 
\Big| \displaystyle{ p_{n-1}  \over  r_{n-1} } -  \displaystyle{ p_{n-3}  \over  r_{n-3} } \Big|  
\over  
| \Im (\alpha^{\prime}_0) | 
} 
< \varepsilon
\]  
 
\noindent
 for all sufficiently large $n$. Therefore, 
$
\Big| \beta_n^{\prime}+\displaystyle{r_{n-2}  \over r_{n-1} }\alpha_n^{\prime}
+\displaystyle{r_{n-3}  \over r_{n-1} } \Big| 
< 
\displaystyle{ r_{n-3}  \over  r_{n-1} } \varepsilon 
< 
\varepsilon  
$   
so our theorem is proved. \hfill \qed

\bigskip

Since   
$\Im( \overline{\alpha_0^{\prime} } )=-\Im(\alpha_0^{\prime})$, where $\overline{\alpha_0^{\prime} }$
is the complex conjugate of $\alpha_0^{\prime} $,   
we can change the sign of  
$\Delta_N \Delta_N^{\prime} \Im(\alpha_0^{\prime})  \Im(\beta_0^{\prime} ) $ 
by changing $\alpha_0^{\prime}$ for $\overline{\alpha_0^{\prime} } $, if necessary,   
and just the same with  $\beta_0$.  

When $\alpha_0^{\prime}$ is real and $\beta_0^{\prime}$ is imaginary,  
we can prove the same result as the above theorem in a similar way 
under the assumption  
$\Delta_N \Delta_N^{\prime} (\alpha_0-\alpha_0^{\prime}) ( \beta_0 - \Re (\beta_0^{\prime} ) )<0$, 
where $\Re(*)$ denotes the real part of $*$.   
Further, when $\alpha_0^{\prime}$ is imaginary and $\beta_0^{\prime}$ is real,  
it can be seen  that  the condition   
$\Delta_N \Delta_N^{\prime} (\alpha_0- \Re (\alpha_0^{\prime}) ) ( \beta_0 - \beta_0^{\prime}  )<0$  
leads to the same result.    
We may consider these results the dimension 2-version of (3.1). 

\vspace{1cm}

{\bf 4.}\  The boundedness of $\beta_n$.     
 
\medskip

Let $\mathcal{D}$ and $\mathcal{D}^{\prime}$ be as before. 
Let $(a_n,\,b_n)\ (0 \leq a_n \leq b_n,\ 1 \leq b_n,\ n \geq 0)$ be a sequence with the admissibility. 
We define rational functions with two variables $X$ and  $Y$
\[ F
\left[
\begin{array}{@{\,}c@{\,}}
a_0    \\
b_0   
\end{array}
\right]  
\left(
\begin{array}{@{\,}c@{\,}}
X    \\
Y   
\end{array}
\right)  
=
\left(
\begin{array}{@{\,}c@{\,}}
a_0+X    \\
b_0+Y   
\end{array}
\right),  
\] 
and for $n \geq 1$     
\[ F
\left[
\begin{array}{@{\,}cccc@{\,}}
a_0  & a_1 & \cdots & a_n   \\
b_0  & b_1 & \cdots & b_n   
\end{array}
\right]  
\left(
\begin{array}{@{\,}c@{\,}}
X    \\
Y   
\end{array}
\right)  
=
 F
\left[
\begin{array}{@{\,}cccc@{\,}}
a_0  & a_1 & \cdots & a_{n-1}   \\
b_0  & b_1 & \cdots & b_{n-1}   
\end{array}
\right]  
\left(
\begin{array}{@{\,}c@{\,}}
1/(b_n+Y)    \\
(a_n+X)/(b_n+Y)   
\end{array}
\right).   
\]
Then, for example 
\[ F
\left[
\begin{array}{@{\,}cc@{\,}}
a_0  & a_1  \\
b_0  & b_1 
\end{array}
\right]  
\left(
\begin{array}{@{\,}c@{\,}}
X    \\
Y   
\end{array}
\right)  
=
\left(
\begin{array}{@{\,}c@{\,}}
a_0+1/(b_1+Y)    \\
b_0+(a_1+X)/(b_1+Y)   
\end{array}
\right).   
\]
When $X=x,\ Y=y$ ($0< x < 1,\ 0 < y <1$ if $a_1 < b_1$, and $0 < x < y < 1$ if $a_1 =b_1$),  
\[
\psi \left(
F
\left[
\begin{array}{@{\,}cc@{\,}}
a_0  & a_1  \\
b_0  & b_1 
\end{array}
\right]  
\left(
\begin{array}{@{\,}c@{\,}}
x    \\
y   
\end{array}
\right)  \right)
=
\left(
\begin{array}{@{\,}c@{\,}}
a_1+x    \\
b_1+y   
\end{array}
\right)  
=
 F
\left[
\begin{array}{@{\,}c@{\,}}
a_1    \\
b_1   
\end{array}
\right]  
\left(
\begin{array}{@{\,}c@{\,}}
x    \\
y   
\end{array}
\right)  
\]
here $\psi$ is the Jacobi-Perron map defined in the introduction.   
Further we have 
\[ \psi \left( F
\left[
\begin{array}{@{\,}cccc@{\,}}
a_0  & a_1 & \cdots & a_n   \\
b_0  & b_1 & \cdots & b_n   
\end{array}
\right]  
\left(
\begin{array}{@{\,}c@{\,}}
x    \\
y   
\end{array}
\right)  \right)
= 
F
\left[
\begin{array}{@{\,}cccc@{\,}}
a_1  & a_2 & \cdots & a_n   \\
b_1  & b_2 & \cdots & b_n   
\end{array}
\right]  
\left(
\begin{array}{@{\,}c@{\,}}
x    \\
y   
\end{array}
\right),  
\]   
\[ {\psi}^n \left( F
\left[
\begin{array}{@{\,}cccc@{\,}}
a_0  & a_1 & \cdots & a_n   \\
b_0  & b_1 & \cdots & b_n   
\end{array}
\right]  
\left(
\begin{array}{@{\,}c@{\,}}
x    \\
y   
\end{array}
\right)  \right)
= 
F
\left[
\begin{array}{@{\,}c@{\,}}
a_n   \\
b_n 
\end{array}
\right]  
\left(
\begin{array}{@{\,}c@{\,}}
x    \\
y   
\end{array}
\right)  
=
\left(
\begin{array}{@{\,}c@{\,}}
a_n+x    \\
b_n+y   
\end{array}
\right). 
\]

For $n \geq 0$, we denote by 
\[
D
\left(
\begin{array}{@{\,}cccc@{\,}}
a_0 & a_{1} & \cdots & a_n   \\
b_0 & b_{1} & \cdots & b_n  
\end{array}
\right)  
\]
the subset of $\mathcal{D}$ 
\[
\left\{
F
\left[
\begin{array}{@{\,}cccc@{\,}}
a_0  & a_1 & \cdots & a_n   \\
b_0  & b_1 & \cdots & b_n   
\end{array}
\right]  
\left(
\begin{array}{@{\,}c@{\,}}
x    \\
y   
\end{array}
\right)  
\ ;\ 
\begin{array}{@{\,}l@{\,}}
0 \leq x \leq 1,   0 \leq y \leq 1  \ {\rm when}\  a_n < b_n \ {\rm and}\ \\
0 \leq x \leq y \leq 1 \ {\rm when}\  a_n = b_n,   \ {\rm respectively}
\end{array}
\right\}.
\]
For $(x,\,y)$ ($0<x<1,\ 0<y<1$ when $a_n<b_n$ and $0<x<y<1$ when $a_n=b_n$, respectively),  
the expansion from the first step to $n+1$th step of the point 
\[
F
\left[
\begin{array}{@{\,}cccc@{\,}}
a_0  & a_1 & \cdots & a_n   \\
b_0  & b_1 & \cdots & b_n   
\end{array}
\right]  
\left(
\begin{array}{@{\,}c@{\,}}
x    \\
y   
\end{array}
\right)  \ \ \ \in\ \ \mathcal{D}
\] 
is 
\[
\left[
\begin{array}{@{\,}cccc@{\,}}
a_0 & a_1 & \cdots & a_n  \\
b_0 & b_1 & \cdots & b_n 
\end{array}
\right].  
\]
The next Lemma can be seen easily by considerimg the union of the quadrangles with vertexes  
 $\big( \displaystyle{ 1  \over  b+1 },\ \displaystyle{ a  \over  b+1 } \big)$,\ 
  $\big( \displaystyle{ 1  \over  b+1 },\ \displaystyle{ a+1  \over  b+1 } \big)$,\ 
  $\big( \displaystyle{ 1  \over  b },\ \displaystyle{ a  \over  b } \big)$,\ 
 $\big( \displaystyle{ 1  \over  b },\ \displaystyle{ a+1  \over  b } \big)$  
when $a<b$, and that of triangles with vertexes  
 $\big( \displaystyle{ 1  \over  b+1 },\ \displaystyle{ b  \over  b+1 } \big)$,\ 
  $\big( \displaystyle{ 1  \over  b+1 },\ 1 \big)$,\ 
  $\big( \displaystyle{ 1  \over  b },\ 1 \big)$ when $a=b$.     

\bigskip

{\bf Lemma 4.1.}\ \  {\it  It holds that }          
\[
\left\{
\left(
\begin{array}{@{\,}c@{\,}}
x   \\
y 
\end{array} 
\right)  
\ ;\   
0 < x \leq 1,\ 0 \leq y \leq 1 
\right\}
\]  
\[=  
\left(
\bigcup_{0 \leq a < b } 
\left\{
\left(
\begin{array}{@{\,}c@{\,}}
1/(b+y)   \\
 (a+x)/(b+y) 
\end{array} 
\right)
\ ;\ 
0 \leq x \leq 1,\ 0 \leq y \leq 1 
\right\} 
\right)
\]  
\[
\bigcup
\left(
\bigcup_{1 \leq a = b }  
\left\{
\left(
\begin{array}{@{\,}c@{\,}}
1/(b+y)   \\
 (a+x)/(b+y) 
\end{array} 
\right)
\ ;\ 
0 \leq x \leq y \leq 1 
\right\}
\right), 
\]  
{\it and} 
\[
\left\{
\left(
\begin{array}{@{\,}c@{\,}}
x   \\
y 
\end{array} 
\right)  
\ ;\   
0 < x \leq  y \leq 1 
\right\}
\]  
\[=  
\left(
\bigcup_{1 \leq a < b } 
\left\{
\left(
\begin{array}{@{\,}c@{\,}}
1/(b+y)   \\
 (a+x)/(b+y) 
\end{array} 
\right)
\ ;\ 
0 \leq x \leq 1,\ 0 \leq y \leq 1 
\right\} 
\right)
\]  
\[
\bigcup
\left(
\bigcup_{1 \leq a = b }  
\left\{
\left(
\begin{array}{@{\,}c@{\,}}
1/(b+y)   \\
 (a+x)/(b+y) 
\end{array} 
\right)
\ ;\ 
0 \leq x \leq y \leq 1 
\right\}
\right), 
\]  
{\it  where a and b run over the rational integers with the conditions respectively. 
}

\bigskip

{\bf Lemma 4.2.}\ \  {\it  Let $|\ \ast \ |$ be the measure of a domain $ \ast $. We have }   

\medskip

$ 
({\rm i})\ \ \left|
D
\left(
\begin{array}{@{\,}cccc@{\,}}
a_0 & a_{1} & \cdots & a_{n-1}   \\
b_0 & b_{1} & \cdots & b_{n-1}  
\end{array}
\right)  \right|
=  
\displaystyle\sum_{0 \leq a \leq b,\ 1 \leq b} 
\left| D
\left(
\begin{array}{@{\,}ccccc@{\,}}
a_0 & a_{1} & \cdots & a_{n-1}  &  a  \\
b_0 & b_{1} & \cdots & b_{n-1}  &  b
\end{array}
\right) \right|
$

\medskip

\ \ \ \ {\it when}\   $a_{n-1} < b_{n-1}$, {\it and} 

\medskip

$ 
({\rm ii})\ \ 
\left| D
\left(
\begin{array}{@{\,}cccc@{\,}}
a_0 & a_{1} & \cdots & a_{n-1}   \\
b_0 & b_{1} & \cdots & b_{n-1}  
\end{array}
\right)  \right|
= 
\displaystyle\sum_{1 \leq a \leq b } 
\left| D
\left(
\begin{array}{@{\,}ccccc@{\,}}
a_0 & a_{1} & \cdots & a_{n-1}  &  a  \\
b_0 & b_{1} & \cdots & b_{n-1}  &  b
\end{array}
\right) \right|
$

\medskip

\ \ \ \ {\it when}\   $a_{n-1} = b_{n-1}$. 

\medskip

{\it Proof.}\ \  
(i)\  
\[
\bigcup_{0 \leq a \leq b,\ 1 \leq b} 
D
\left(
\begin{array}{@{\,}ccccc@{\,}}
a_0 & a_{1} & \cdots & a_{n-1}  &  a  \\
b_0 & b_{1} & \cdots & b_{n-1}  &  b
\end{array}
\right)  
\] 
\[ =                           
\left(
\bigcup_{0 \leq a < b } 
D
\left(
\begin{array}{@{\,}ccccc@{\,}}
a_0 & a_{1} & \cdots & a_{n-1}  &  a  \\
b_0 & b_{1} & \cdots & b_{n-1}  &  b
\end{array}
\right)  
\right) 
\bigcup   
\left(
\bigcup_{1 \leq a = b } 
D
\left(
\begin{array}{@{\,}ccccc@{\,}}
a_0 & a_{1} & \cdots & a_{n-1}  &  a  \\
b_0 & b_{1} & \cdots & b_{n-1}  &  b
\end{array}
\right)  
\right) 
\]  
\[= 
\left(
\bigcup_{0 \leq a < b}
\left\{
F
\left[
\begin{array}{@{\,}cccc@{\,}}
a_0  & a_1 & \cdots & a_{n-1}    \\
b_0  & b_1 & \cdots & b_{n-1}     
\end{array}
\right]  
\left(
\begin{array}{@{\,}c@{\,}}
1/(b+y)    \\
(a+x)/(b+y)   
\end{array}
\right)  
\ ;\ 
0 \leq x \leq 1,   0 \leq y \leq 1  
\right\}  
\right) 
\]  
\[
\bigcup 
\left(
\bigcup_{1 \leq a = b}
\left\{
F
\left[
\begin{array}{@{\,}cccc@{\,}}
a_0  & a_1 & \cdots & a_{n-1}    \\
b_0  & b_1 & \cdots & b_{n-1}    
\end{array}
\right]  
\left(
\begin{array}{@{\,}c@{\,}}
1/(b+y)    \\
(a+x)/(b+y)   
\end{array}
\right)  
\ ;\ 
0 \leq x \leq y \leq 1  
\right\}  
\right)
\]  
\[
= 
\left\{
F
\left[
\begin{array}{@{\,}cccc@{\,}}
a_0  & a_1 & \cdots & a_{n-1}    \\
b_0  & b_1 & \cdots & b_{n-1}    
\end{array}
\right]  
\left(
\begin{array}{@{\,}c@{\,}}
x    \\
y   
\end{array}
\right)  
\ ;\ 
0 < x \leq 1, \  0 \leq y \leq 1  
\right\}\ \ {\rm by\ Lemma\ 4.1.}  
\]  
The measure of this set  equals the left hand side of (i).   
Note that  for $a$ and $b$ 
\[ 
D
\left(
\begin{array}{@{\,}ccccc@{\,}}
a_0 & a_{1} & \cdots & a_{n-1}   & a  \\
b_0 & b_{1} & \cdots & b_{n-1}   & b  
\end{array}
\right)  \ \  
(0 \leq a \leq b,\ 1 \leq b)
\]  
are mutually disjoint except that they may have only boundaries in common. 

(ii)
\  
\[
\bigcup_{1 \leq a \leq b } 
D
\left(
\begin{array}{@{\,}ccccc@{\,}}
a_0 & a_{1} & \cdots & a_{n-1}  &  a  \\
b_0 & b_{1} & \cdots & b_{n-1}  &  b
\end{array}
\right)  
\] 
\[ =                           
\left(
\bigcup_{1 \leq a < b } 
D
\left(
\begin{array}{@{\,}ccccc@{\,}}
a_0 & a_{1} & \cdots & a_{n-1}  &  a  \\
b_0 & b_{1} & \cdots & b_{n-1}  &  b
\end{array}
\right)  
\right) 
\bigcup   
\left(
\bigcup_{1 \leq a = b } 
D
\left(
\begin{array}{@{\,}ccccc@{\,}}
a_0 & a_{1} & \cdots & a_{n-1}  &  a  \\
b_0 & b_{1} & \cdots & b_{n-1}  &  b
\end{array}
\right)  
\right) 
\]  
\[= 
\left(
\bigcup_{1 \leq a < b}
\left\{
F
\left[
\begin{array}{@{\,}cccc@{\,}}
a_0  & a_1 & \cdots & a_{n-1}    \\
b_0  & b_1 & \cdots & b_{n-1}     
\end{array}
\right]  
\left(
\begin{array}{@{\,}c@{\,}}
1/(b+y)    \\
(a+x)/(b+y)   
\end{array}
\right)  
\ ;\ 
0 \leq x \leq 1,   0 \leq y \leq 1  
\right\}  
\right) 
\]  
\[
\bigcup 
\left(
\bigcup_{1 \leq a = b}
\left\{
F
\left[
\begin{array}{@{\,}cccc@{\,}}
a_0  & a_1 & \cdots & a_{n-1}    \\
b_0  & b_1 & \cdots & b_{n-1}    
\end{array}
\right]  
\left(
\begin{array}{@{\,}c@{\,}}
1/(b+y)    \\
(a+x)/(b+y)   
\end{array}
\right)  
\ ;\ 
0 \leq x \leq y \leq 1  
\right\}  
\right)
\]  

$
= 
\left\{
F
\left[
\begin{array}{@{\,}cccc@{\,}}
a_0  & a_1 & \cdots & a_{n-1}    \\
b_0  & b_1 & \cdots & b_{n-1}    
\end{array}
\right]  
\left(
\begin{array}{@{\,}c@{\,}}
x    \\
y   
\end{array}
\right)  
\ ;\ 
0 < x \leq  y \leq 1  
\right\}\ \ {\rm by\ Lemma\ 4.1.}  
$ 

\medskip
\noindent  
The measure of this set is equals the left hand side of (ii).   
\hfill  \qed

\bigskip

We see that 
$ 
\psi\ :\ D^{\circ} \bigl(^a_b \bigr)\  \longrightarrow\ \mathcal{D}^{\circ}  
( \ {\rm or}\  (\mathcal{D}^{\prime})^{\circ} \ {\rm when}\  a=b), 
$ 
where $D^{\circ} \bigl(^a_b \bigr)$ is the  set of all interior points of $D \bigl(^a_b \bigr)$, 
is a bijection and   also  $\psi$ maps
$
D^{\circ}
\left(
\begin{array}{@{\,}cccc@{\,}}
a_0 & a_{1} & \cdots & a_{n}   \\
b_0 & b_{1} & \cdots & b_{n}  
\end{array}
\right)
$ 
onto 
$
D^{\circ}
\left(
\begin{array}{@{\,}cccc@{\,}}
a_1 & a_{2} & \cdots & a_{n}   \\
b_1 & b_{2} & \cdots & b_{n}  
\end{array}
\right)
$ 
bijectively. 

\bigskip

{\bf Lemma 4.3.}\ \  {\it We consider  the bijection   
$\psi \ :\ D^{\circ} \bigl(^a_b \bigr)  \longrightarrow  \mathcal{D}^{\circ} $
( or $(\mathcal{D}^{\prime})^{\circ}$  when  $a=b$). 
The image of a segment in $D^{\circ} \bigl(^a_b \bigr)$ is a segment in 
$\mathcal{D}^{\circ} $( or $ (\mathcal{D}^{\prime})^{\circ}$) and  
the inverse image of a segment in $\mathcal{D}^{\circ} $( or $(\mathcal{D}^{\prime})^{\circ} $)   
is a segment in $D^{\circ} \bigl(^a_b \bigr)$.  
The inverse image of a convex set in $\mathcal{D}^{\circ} $( or $ (\mathcal{D}^{\prime})^{\circ}$)   
is convex. 
}   

\medskip

{\it Proof.}\ \ 
Let $P_1,\ P_2\ (P_1 \not= P_2)$ be any two points in $D^{\circ} \bigl(^a_b \bigr)$, 
and let $Q_1=\psi(P_1)=(x_1,\,y_1)$, $Q_2=\psi(P_2)=(x_2,\,y_2)$.  
Note that $0<x_i < y_i$, $1<y_i\ ( i=1,2)$. 
It holds 
$P_1=(a+1/{y_1},\ b+{x_1}/{y_1})$ and $P_2=(a+1/{y_2},\ b+{x_2}/{y_2})$. 
Let $P$ be any point on the segment $P_1 P_2$, then   
\[
P=
\Big(
s(a+1/{y_1})+(1-s)(a+1/{y_2}),\ 
s(b+{x_1}/{y_1})+(1-s)(b+{x_2}/{y_2})  
\Big) 
\] 
by means of  a parameter $s\ (0 \leq s \leq 1)$. 
Let 
$t=sy_2/((1-s)y_1+sy_2)$,  
then $0 \leq t \leq 1$ and $s=ty_1/(ty_1+(1-t)y_2)$. 
From this, we have 
\[
P=
\Big(
a+ 1/(ty_1+(1-t)y_2),\ 
b+ (tx_1+(1-t)x_2)/(ty_1+(1-t)y_2)
\Big)  
\]
so, 
\[
\psi(P)=
(
tx_1+(1-t)x_2,\  
ty_1+(1-t)y_2
).  
\]
Therefore, $\psi(P)$ lies on the segment $Q_1 Q_2$. 
Because the parameter $t$ runs from 0 to 1 when $s$ runs from 0 to 1, 
the segment $P_1 P_2$ is mapped onto $Q_1 Q_2$. 
Since $\psi$ is bijection, the inverse is clearly the case.  

Let $S$ be a convex set in $\mathcal{D}^{\circ} $( or $(\mathcal{D}^{\prime})^{\circ} $).   
Let $P_1$ and $P_2$ $(P_1 \not= P_2)$ be any two points in $\psi^{-1}(S)$, 
the inverse image of $S$, 
and let 
$Q_1=\psi(P_1)$,  $Q_2=\psi(P_2)$. 
The image of $P$, any point on the segment $P_1 P_2$, is on the segment $Q_1 Q_2$, 
so it belongs to $S$. 
Thus $P \in \psi^{-1}(S)$, and $\psi^{-1}(S)$ is proved to be convex. 
  \hfill  \qed
 
\bigskip 

Note that  this lemma is the case for 
\[   
\psi\ :\ 
D^{\circ}
\left(
\begin{array}{@{\,}cccc@{\,}}
a_0 & a_{1} & \cdots & a_{n}   \\
b_0 & b_{1} & \cdots & b_{n}  
\end{array}
\right) 
\ \longrightarrow  \ 
D^{\circ}
\left(
\begin{array}{@{\,}cccc@{\,}}
a_1 & a_{2} & \cdots & a_{n}   \\
b_1 & b_{2} & \cdots & b_{n}  
\end{array}
\right)  
\] 
(those are subsets of $D^{\circ} \bigl(^{a_0}_{b_0} \bigr)$ and $\mathcal{D}^{\circ} $
( or $(\mathcal{D}^{\prime})^{\circ}$) ),  respectively,  
and means that the boundary of  any 
\[ 
D
\left(
\begin{array}{@{\,}cccc@{\,}}
a_0 & a_{1} & \cdots & a_{n}   \\
b_0 & b_{1} & \cdots & b_{n}  
\end{array}
\right)  
\]  
is a quadrangle if $a_n < b_n$ and a triangle if $a_n = b_n$, respectively. 
Its vertexes are the points according to 
$(X,\,Y)=(0,0),\ (1,0),\ (1,1),\ (0,1)$
  and $(X,\,Y)=(0,0),\ (1,1),\ (0,1)$ 
of 
\[
F
\left[
\begin{array}{@{\,}cccc@{\,}}
a_0  & a_1 & \cdots & a_{n}    \\
b_0  & b_1 & \cdots & b_{n}    
\end{array}
\right]  
\left(
\begin{array}{@{\,}c@{\,}}
X    \\
Y   
\end{array} 
\right),  
\]  
respectively. 
We denote this function by $P_n(X,Y)$ (depends on the admissible sequence $(a_i\,b_i)\ ( 0 \leq i \leq n)$)  
and denote vertexes by $P_n(0,0),\ P_n(1,0),\ P_n(1,1)$ and $ P_n(0,1)$ accordingly.   

Let $a,\ b,\ c,\ d$ ($0 \leq a \leq b,\ 1 \leq b$, $0 \leq c \leq d,\ 1 \leq d$) be integers, 
where $1 \leq c$ if $a=b$ (admissibility).   
Let $P_i=\bigl(a+1/(d+y_i),\ b+(c+x_i)/(d+y_i) \bigr)$ $(i=1,2,3)$ where 
$0<x_i <1,\ 0<y_i <1$ if $c<d$, and $0< x_i < y_i <1$ if $c=d$.   
Then, $P_i \in D^{\circ} \bigl(^a_b \bigr)$. 
Let $Q_i=\bigl(c+x_i,\ d+y_i \bigr)$, 
then $Q_i \in D^{\circ} \bigl(^c_d \bigr)$ and $\psi (P_i)=Q_i$.  
The triangle $\triangle P_1 P_2 P_3$   
and its interior points  are mapped bijectively by $\psi$  on the triangle $\triangle Q_1 Q_2 Q_3$  and its interior points.   
Note that, in general,   
$D^{\circ} \bigl(^{a\ c}_{b\ d} \bigr)=D^{\circ} \bigl(^a_b \bigr) \cap \psi^{-1}(D^{\circ} \bigl(^c_d \bigr))$.  
We denote by $|\triangle P_1 P_2 P_3|$ the measure of $\triangle P_1 P_2 P_3$ as before.    
Now, we get the following. 

\bigskip

{\bf Lemma 4.4.}\ \  
{\it $|\triangle P_1 P_2 P_3|   
= \displaystyle{ 1  \over  (d+y_1) (d+y_2) (d+y_3) } |\triangle Q_1 Q_2 Q_3|.$ }

\medskip

{\it Proof.}\ \ It holds 
\[|\triangle Q_1 Q_2 Q_3|=\displaystyle{1  \over  2} 
{\rm abs} 
\left|          
   \begin{array}{ccc}
     1  &  1  &  1   \\
     c+x_1  & c+x_2   & c+x_3     \\
     d+y_1  &  d+y_2  &  d+y_3   
    \end{array}
   \right| 
=  
\displaystyle{1  \over  2} 
{\rm abs} 
\left|          
   \begin{array}{ccc}
     1  &  1  &  1   \\
     x_1  & x_2   & x_3     \\
     y_1  &  y_2  &  y_3   
    \end{array}
   \right|,     
\]    
and  

\medskip

$
|\triangle P_1 P_2 P_3| 
=  
\displaystyle{1  \over  2} 
{\rm abs} 
\left|          
   \begin{array}{ccc}
     1  &  1  &  1   \\  [4pt]
    a+\displaystyle{ 1  \over  d+y_1}  & a+\displaystyle{ 1  \over  d+y_2}   & a+\displaystyle{ 1  \over  d+y_3}     \\  [12pt]
    b+\displaystyle{ c+x_1  \over  d+y_1}  &   b+\displaystyle{ c+x_2  \over  d+y_2}   &  b+\displaystyle{ c+x_3  \over  d+y_3}     
    \end{array}
   \right|     
$
   
$=
\displaystyle{1  \over  2} 
{\rm abs} 
\left|          
   \begin{array}{ccc}
     1  &  1  &  1   \\  [4pt]
    \displaystyle{ 1  \over  d+y_1}  & \displaystyle{ 1  \over  d+y_2}   & \displaystyle{ 1  \over  d+y_3}     \\   [12pt]
    \displaystyle{ c+x_1  \over  d+y_1}  &   \displaystyle{ c+x_2  \over  d+y_2}   &  \displaystyle{ c+x_3  \over  d+y_3}     
    \end{array}
   \right|  
$  
  
$
= 
\displaystyle{1  \over  2} 
\displaystyle{ 1  \over  (d+y_1) (d+y_2) (d+y_3) }
{\rm abs} 
\left|          
   \begin{array}{ccc}
    d+y_1  &  d+y_2  &  d+y_3   \\
    1   &   1  &  1    \\
    c+x_1   &    c+x_2    &   c+x_3      
    \end{array}
   \right|  
$    

$
=\displaystyle{ 1  \over  (d+y_1) (d+y_2) (d+y_3) } |\triangle Q_1 Q_2 Q_3|. 
$  
 \hfill  \qed

\bigskip

{\bf Lemma 4.5.}\ \  {\it  Let $a,\,b,\,c,$ and $d$ be as above.      
Suppose a domain in $D^{\circ} \bigl(^c_d \bigr)$ has the measure $N>0$ and  
its inverse image by $\psi$ in $D^{\circ} \bigl(^a_b \bigr)$ has the measure $N^{\prime}>0$,  then 
\[ 
\displaystyle{  1  \over  (d+1)^3 } N <  N^{\prime}   
< \displaystyle{  1  \over  d^3 } N. 
\]
}  
\medskip
{\it Proof.}\ \ In the case of Lemma 4.4,  we see 
\[ 
 \displaystyle{ 1  \over (d+1)^3 } |\triangle Q_1 Q_2 Q_3| < |\triangle P_1 P_2 P_3|  
 < \displaystyle{ 1  \over d^3 } |\triangle Q_1 Q_2 Q_3|.
\]   
Because $\triangle Q_1 Q_2 Q_3$ is any triangle in $D^{\circ} \bigl(^c_d \bigr)$, the lemma follows. 
 \hfill  \qed 

\bigskip

Let $\{ r_n \}$,\ $\{ p_n \}$ and $\{ q_n \}\ (n \geq 0)$ be sequences defined by (1.2) 
according to $(a_n,\,b_n)$ $ (n \geq 0)$.  

\bigskip

{\bf  Lemma 4.6.}\ \  {\it It holds for $n \geq 0$}   
\[ 
F
\left[
\begin{array}{@{\,}cccc@{\,}}
a_0  & a_1 & \cdots & a_{n}    \\
b_0  & b_1 & \cdots & b_{n}    
\end{array}
\right]  
\left(
\begin{array}{@{\,}c@{\,}}
X    \\
Y   
\end{array} 
\right)  
= 
\left( 
\displaystyle{p_{n-2}X+p_{n-1}Y+p_n   \over  r_{n-2}X+r_{n-1}Y+r_n  },\ 
\displaystyle{q_{n-2}X+q_{n-1}Y+q_n   \over  r_{n-2}X+r_{n-1}Y+r_n  } 
\right).  
\]  

\medskip

{\it Proof.}\ \ Since, from Introduction, 
\[ K_0=\left(          
   \begin{array}{ccc}
     0  &   0  &  1   \\
     1  &   0  &  a_0  \\
     0  &   1  &  b_0   
    \end{array}
   \right)
=\left(          
   \begin{array}{ccc}
     r_{-2}  &   r_{-1}  &  r_0   \\
     p_{-2}  &   p_{-1}  &  p_0   \\
     q_{-2}  &   q_{-1}  &  q_0  
    \end{array}
   \right), 
\]
we see that our theorem holds for $n=0$.  
Assume  the case of $n \geq 0$ is true, then 
\[ 
F
\left[
\begin{array}{@{\,}cccc@{\,}}
a_0  & a_1 & \cdots & a_{n+1}    \\
b_0  & b_1 & \cdots & b_{n+1}    
\end{array}
\right]  
\left(
\begin{array}{@{\,}c@{\,}}
X    \\
Y   
\end{array} 
\right)  
= 
F
\left[
\begin{array}{@{\,}cccc@{\,}}
a_0  & a_1 & \cdots & a_{n}    \\
b_0  & b_1 & \cdots & b_{n}    
\end{array}
\right]  
\left(
\begin{array}{@{\,}c@{\,}}
1/(b_{n+1}+Y)    \\
(a_{n+1}+X)/(b_{n+1}+Y)   
\end{array} 
\right)  
\] 
\[ = 
\left( 
\begin{array}{@{\,}c@{\,}}
\displaystyle{p_{n-2}/(b_{n+1}+Y)+p_{n-1} (a_{n+1}+X)/(b_{n+1}+Y)+p_n   \over  r_{n-2}/(b_{n+1}+Y)+r_{n-1} (a_{n+1}+X)/(b_{n+1}+Y)+r_n  }   \\  [12pt]
\displaystyle{ q_{n-2}/(b_{n+1}+Y)+q_{n-1} (a_{n+1}+X)/(b_{n+1}+Y)+q_n   \over  r_{n-2}/(b_{n+1}+Y)+r_{n-1}( a_{n+1}+X)/(b_{n+1}+Y)+r_n  } 
\end{array} 
\right)  
\]  
\[ = 
\left( 
\begin{array}{@{\,}c@{\,}}
\displaystyle{p_{n-2}+p_{n-1} (a_{n+1}+X)+p_n(b_{n+1}+Y)   \over  r_{n-2}+r_{n-1} (a_{n+1}+X)+r_n (b_{n+1}+Y) }   \\   [12pt]
\displaystyle{ q_{n-2}+q_{n-1} (a_{n+1}+X)+q_n (b_{n+1}+Y)  \over  r_{n-2}+r_{n-1}( a_{n+1}+X)+r_n (b_{n+1}+Y) } 
\end{array} 
\right)  
\]  
\[ = 
\left( 
\begin{array}{@{\,}c@{\,}}
\displaystyle{p_{n-1}X+p_nY+p_{n-2}+p_{n-1} a_{n+1}+p_n b_{n+1}   \over  r_{n-1}X+r_nY+r_{n-2}+r_{n-1} a_{n+1}+r_n b_{n+1} }   \\   [12pt]
\displaystyle{ q_{n-1}X+q_nY+q_{n-2}+q_{n-1} a_{n+1}+q_n b_{n+1}   \over  r_{n-1}X+r_nY+r_{n-2}+r_{n-1} a_{n+1}+r_n b_{n+1} } 
\end{array} 
\right)  
\]  
\[ = 
\left( 
\displaystyle{p_{n-1}X+p_nY+p_{n+1}   \over  r_{n-1}X+r_nY+r_{n+1}  }, \   
\displaystyle{ q_{n-1}X+q_nY+q_{n+1}   \over  r_{n-1}X+r_nY+r_{n+1} } 
\right).   
\]  
By induction, our theorem is proved.  
 \hfill  \qed

\bigskip

{\bf Lemma 4.7.}\ \  {\it  Let  $(a_n,\,b_n)\ (n \geq 0)$    
 be a sequence with the admissibility and $P_n$ the function as before.  
Points $P_n(t,0)$, $P_n(t,1)$ and $P_n(t,t)\ (0 \leq t \leq 1)$ exist on the segments 
$P_n(0,0)P_n(1,0)$, $P_n(0,1)P_n(1,1)$ and $P_n(0,0)P_n(1,1)$, respectively. }

\medskip

{\it Proof.}\ \ Since  
\[ 
P_n(0,0)=\left(\displaystyle{p_n  \over  r_n},\ \displaystyle{q_n  \over  r_n} \right),\ \ \ 
P_n(t,0)=\left(\displaystyle{p_{n-2}t+p_n  \over  r_{n-2}t+r_n},\ \displaystyle{q_{n-2}t+q_n  \over  r_{n-2}t+r_n} \right),  
\] 
the slope of the segment $P_n(0,0)P_n(t,0)$ is 
\[ 
\bigg( \displaystyle{q_{n-2}t+q_n  \over  r_{n-2}t+r_n} - \displaystyle{q_n  \over  r_n} \bigg) 
\bigg/ 
\bigg( \displaystyle{p_{n-2}t+p_n  \over  r_{n-2}t+r_n} - \displaystyle{p_n  \over  r_n} \bigg) 
= 
\displaystyle{ r_n q_{n-2}-r_{n-2} q_n   \over   r_n p_{n-2}-r_{n-2} p_n }. 
\] 
This is independent of $t$, so $P_n(t,0)$ exists on the segment $P_n(0,0)P_n(1,0)$. 

\noindent 
The slope of the segment $P_n(0,1)P_n(t,1)$ is 
\[ 
\bigg( \displaystyle{q_{n-2}t+q_{n-1}+q_n  \over  r_{n-2}t+r_{n-1}+r_n} - \displaystyle{ q_{n-1}+q_n  \over  r_{n-1}+r_n } \bigg)
\bigg/ 
\bigg( \displaystyle{p_{n-2}t+p_{n-1}+p_n  \over  r_{n-2}t+r_{n-1}+r_n} - \displaystyle{p_{n-1}+p_n  \over r_{n-1}+ r_n} \bigg) 
\]  
\[
= 
\displaystyle{ (r_{n-2} q_{n-1}-r_{n-1} q_{n-2})+(r_{n-2} q_{n}-r_{n} q_{n-2})   
\over   (r_{n-2} p_{n-1}-r_{n-1} p_{n-2})+ (r_{n-2} p_{n}-r_{n} p_{n-2})}. 
\] 

\noindent
The slope of the segment $P_n(0,0)P_n(t,t)$ is 
\[ 
\bigg( \displaystyle{q_{n-2}t+q_{n-1}t+q_n  \over  r_{n-2}t+r_{n-1}t+r_n} - \displaystyle{ q_n  \over  r_n } \bigg) 
\bigg/ 
\bigg( \displaystyle{p_{n-2}t+p_{n-1}t+p_n  \over  r_{n-2}t+r_{n-1}t+r_n} - \displaystyle{p_n  \over  r_n} \bigg) 
\] 
\[
= 
\displaystyle{ (r_{n-2} q_{n}-r_{n} q_{n-2})+(r_{n-1} q_{n}-r_{n} q_{n-1})   
\over   (r_{n-2} p_{n}-r_{n} p_{n-2})+ (r_{n-1} p_{n}-r_{n} p_{n-1})}. 
\] 
These are also independent of $t$ and our Lemma follows. 
 \hfill  \qed

\bigskip

We denote by $S_n (t)\ (0 \leq t \leq 1)$,  
depending on a natural number $n$ and a sequence $(a_i,\,b_i)\ ( 0 \leq i \leq n )$ with the admissibility, 
 the measure of the quadrangle with vertexes 
$P_n(0,0)$, $P_n(t,0)$, $P_n(t,1)$ and $P_n(0,1)$ 
when $a_n < b_n$, 
and by $S_n^{\prime}(t)\ (0 \leq t \leq 1)$ the measure of the quadrangle with vertexes 
$P_n(0,0)$, $P_n(t,t)$, $P_n(t,1)$ and $P_n(0,1)$ 
when $a_n = b_n$ (note that this is a triangle when $t=1$). 

From the definitions,  
\[ 
S_n (1)=
\left|D
\left(
\begin{array}{@{\,}cccc@{\,}}
a_0 & a_{1} & \cdots & a_n   \\
b_0 & b_{1} & \cdots & b_n  
\end{array}
\right) \right|\ ( a_n < b_n),  
\] 
\[  
S_n^{\prime}(1)=
\left|D
\left(
\begin{array}{@{\,}cccc@{\,}}
a_0 & a_{1} & \cdots & a_n   \\
b_0 & b_{1} & \cdots & b_n  
\end{array}
\right) \right|\ ( a_n = b_n). 
\]

\bigskip

{\bf Theorem 4.1.}\ \  {\it  For  $t\ (0 \leq t \leq 1)$,  we have }    
\[ 
\begin{array}{@{\,}l@{\,}}
S_n(t)=\displaystyle{ t(r_{n-2}t+r_{n-1}+2r_n)   \over   2r_n(r_{n-1}+r_n)(r_{n-2}t+r_n)(r_{n-2}t+r_{n-1}+r_n) },     \\  [12pt]
S_n^{\prime}(t)=\displaystyle{ t(r_{n-2}t+r_{n-1}+r_n (2-t))   \over   2r_n(r_{n-1}+r_n)(r_{n-2}t+r_{n-1}+r_n) (r_{n-2}t+r_{n-1}t+r_n)  }.
\end{array}
\]  

\medskip

{\it Proof.}\ \ The measure of 
$\bigtriangleup  P_n(0,0)P_n(t,0)P_n(t,1)$ is  
\[ 
\displaystyle{ 1  \over  2} 
{\rm abs} 
\left|
\begin{array}{@{\,}ccc@{\,}}
1 & 1 & 1   \\  [8pt]
\displaystyle{p_n  \over  r_n}  &  \displaystyle{p_{n-2}t+p_n  \over  r_{n-2}t+r_n}     &     
   \displaystyle{p_{n-2}t+p_{n-1}+p_n  \over  r_{n-2}t+r_{n-1}+r_n}                \\    [12pt]
\displaystyle{q_n  \over  r_n}  & \displaystyle{q_{n-2}t+q_n  \over  r_{n-2}t+r_n}    &   
   \displaystyle{q_{n-2}t+q_{n-1}+q_n  \over  r_{n-2}t+r_{n-1}+r_n}    
\end{array}
\right|
\]  
\[ \begin{array}{@{\,}l@{\,}}
= \displaystyle{ 1  \over  2r_n (r_{n-2}t+r_n)(r_{n-2}t+r_{n-1}+r_n)} 
{\rm abs} 
\left|
\begin{array}{@{\,}ccc@{\,}}
r_n & r_{n-2}t+r_n &   r_{n-2}t+r_{n-1}+r_n  \\
p_n   &  p_{n-2}t+p_n   &  p_{n-2}t+p_{n-1}+p_n                \\  
q_n   & q_{n-2}t+q_n    &  q_{n-2}t+q_{n-1}+q_n      
\end{array}
\right|    \\
=
\displaystyle{ t  \over  2r_n (r_{n-2}t+r_n)(r_{n-2}t+r_{n-1}+r_n)} 
{\rm abs} 
\left|
\begin{array}{@{\,}ccc@{\,}}
r_{n-2} & r_{n-1} &   r_n  \\
p_{n-2}   &  p_{n-1}  &  p_n                \\  
q_{n-2}   & q_{n-1}   &  q_n      
\end{array}
\right|    \\
=
\displaystyle{ t  \over  2r_n (r_{n-2}t+r_n)(r_{n-2}t+r_{n-1}+r_n)}.  
\end{array}
\] 

The measure of 
$\bigtriangleup  P_n(0,0)P_n(0,1)P_n(t,1)$ is  
\[ 
\displaystyle{ 1  \over  2} 
{\rm abs} 
\left|
\begin{array}{@{\,}ccc@{\,}}
1 & 1 & 1   \\  [8pt]
\displaystyle{p_n  \over  r_n}  &  \displaystyle{p_{n-1}+p_n  \over  r_{n-1}+r_n}     &     
   \displaystyle{p_{n-2}t+p_{n-1}+p_n  \over  r_{n-2}t+r_{n-1}+r_n}                \\  [12pt]  
\displaystyle{q_n  \over  r_n}  & \displaystyle{q_{n-1}+q_n  \over  r_{n-1}+r_n}    &   
   \displaystyle{q_{n-2}t+q_{n-1}+q_n  \over  r_{n-2}t+r_{n-1}+r_n}    
\end{array}
\right|
\]  
\[ \begin{array}{@{\,}l@{\,}}
= \displaystyle{ 1  \over  2r_n (r_{n-1}+r_n)(r_{n-2}t+r_{n-1}+r_n)} 
{\rm abs} 
\left|
\begin{array}{@{\,}ccc@{\,}}
r_n & r_{n-1}+r_n &   r_{n-2}t+r_{n-1}+r_n  \\
p_n   &  p_{n-1}+p_n   &  p_{n-2}t+p_{n-1}+p_n                \\  
q_n   & q_{n-1}+q_n    &  q_{n-2}t+q_{n-1}+q_n      
\end{array}
\right|      \\
 =
\displaystyle{ t  \over  2r_n (r_{n-1}+r_n)(r_{n-2}t+r_{n-1}+r_n)} 
{\rm abs} 
\left|
\begin{array}{@{\,}ccc@{\,}}
r_{n-2} & r_{n-1} &   r_n  \\
p_{n-2}   &  p_{n-1}  &  p_n                \\  
q_{n-2}   & q_{n-1}   &  q_n      
\end{array}
\right|  \\
 =
\displaystyle{ t  \over  2r_n (r_{n-1}+r_n)(r_{n-2}t+r_{n-1}+r_n)}.  
\end{array}
\] 
Since 
\[ 
S_n(t)=|\bigtriangleup  P_n(0,0)P_n(t,0)P_n(t,1)|+|\bigtriangleup  P_n(0,0)P_n(0,1)P_n(t,1) |,
\] 
we get  $S_n(t)$.  
The measure of 
$\bigtriangleup  P_n(0,0)P_n(t,1)P_n(t,t)$ is  
\[ 
\displaystyle{ 1  \over  2} 
{\rm abs} 
\left|
\begin{array}{@{\,}ccc@{\,}}
1 & 1 & 1   \\  [8pt]
\displaystyle{p_n  \over  r_n}  &  \displaystyle{p_{n-2}t+p_{n-1}+p_n  \over  r_{n-2}t+r_{n-1}+r_n }     &     
   \displaystyle{p_{n-2}t+p_{n-1}t+p_n  \over  r_{n-2}t+r_{n-1}t+r_n}                \\    [12pt]
\displaystyle{q_n  \over  r_n}  & \displaystyle{q_{n-2}t+q_{n-1}+q_n  \over  r_{n-2}t+r_{n-1}+r_n }    &   
   \displaystyle{q_{n-2}t+q_{n-1}t+q_n  \over  r_{n-2}t+r_{n-1}t+r_n}    
\end{array}
\right|
\]  
\[ \begin{array}{@{\,}l@{\,}}
=
\displaystyle{ 1  \over  2r_n (r_{n-2}t+r_{n-1}+r_n)(r_{n-2}t+r_{n-1}t+r_n) }    \\
\hspace{2cm}
\times 
{\rm abs} 
\left|
\begin{array}{@{\,}ccc@{\,}}
r_n & r_{n-2}t+r_{n-1}+r_n &   r_{n-2}t+r_{n-1}t+r_n  \\
p_n   &  p_{n-2}t+p_{n-1}+p_n   &  p_{n-2}t+p_{n-1}t+p_n                \\  
q_n   & q_{n-2}t+q_{n-1}+q_n    &  q_{n-2}t+q_{n-1}t+q_n      
\end{array}
\right|    \\
 =
\displaystyle{ t(1-t)  \over  2r_n (r_{n-2}t+r_{n-1}+r_n)(r_{n-2}t+r_{n-1}t+r_n) } 
{\rm abs} 
\left|
\begin{array}{@{\,}ccc@{\,}}
r_{n-2} & r_{n-1} &   r_n  \\
p_{n-2}   &  p_{n-1}  &  p_n                \\  
q_{n-2}   & q_{n-1}   &  q_n      
\end{array}
\right|     \\
 =
\displaystyle{ t(1-t)  \over  2r_n (r_{n-2}t+r_{n-1}+r_n)(r_{n-2}t+r_{n-1}t+r_n)}.  
\end{array}
\] 
Because 
\[ 
S_n^{\prime}(t)=|\bigtriangleup  P_n(0,0)P_n(0,1)P_n(t,1)|+|\bigtriangleup  P_n(0,0)P_n(t,1)P_n(t,t) |,
\] 
we get $S_n^{\prime}(t)$. 
\hfill  \qed

\bigskip

From this theorem we have (for $n \geq 0$) 

\medskip
$  
\left|D
\left(
\begin{array}{@{\,}ccc@{\,}}
a_0 &  \cdots & a_n   \\
b_0 &  \cdots & b_n  
\end{array}
\right) \right| 
= 
S_n (1)
=
\displaystyle{ r_{n-2}+r_{n-1}+2r_n   \over   2r_n(r_{n-1}+r_n)(r_{n-2}+r_n)(r_{n-2}+r_{n-1}+r_n) }    \hfill (4.1)
$  

\medskip
\noindent
when $a_n <b_n$, and  

\medskip

$                  
\left|D
\left(
\begin{array}{@{\,}ccc@{\,}}
a_0 &  \cdots & a_n   \\
b_0 &  \cdots & b_n  
\end{array}
\right) \right| 
= 
S_n^{\prime}(1)
=
\displaystyle{ 1   \over   2r_n(r_{n-1}+r_n)(r_{n-2}+r_{n-1}+r_n)   }     \hfill (4.2)
$ 

\medskip
\noindent
when $a_n = b_n$. 

Next we define the followings ;  

$\mathcal{D}_t \ ( 0 \leq t \leq 1) = \{ (x,\,y) ; t \leq x \leq 1,\ 0 \leq y \leq 1 \}$, 

\medskip
$\mathcal{D}^{\prime}_t \ ( 0 \leq t \leq 1) = \{ (x,\,y) ; t \leq x \leq y \leq 1 \}$, 

\medskip
$(\mathcal{D}_t)^c \ ( 0 \leq t \leq 1) = \{ (x,\,y) ; 0 \leq x \leq t,\ 0 \leq y \leq 1 \}$, 

\medskip
$(\mathcal{D}^{\prime}_t)^c \ ( 0 \leq t \leq 1) = \{ (x,\,y) ; 0 \leq x \leq t,\   x \leq y \leq 1 \}$.  

\medskip
\noindent
Note that the measure of the image of $(\mathcal{D}_t)^c$ by the function 
\[ F
\left[
\begin{array}{@{\,}cccc@{\,}}
a_0  & a_1 & \cdots & a_n   \\
b_0  & b_1 & \cdots & b_n   
\end{array}
\right]  
\left(
\begin{array}{@{\,}c@{\,}}
X    \\
Y   
\end{array}
\right)  
\] 
is $S_n(t)$ and that  of $(\mathcal{D}^{\prime}_t)^c$ is $S_n^{\prime}(t)$.  

If $(\alpha_0,\,\beta_0) \in \mathcal{D}\ (\alpha_0,\,\beta_0$ are irrational numbers)  has  the 
Jacobi-Perron expansion  

$ 
\left[
\begin{array}{@{\,}ccccc@{\,}}
a_0  & a_1 & \cdots & a_n  &  \cdots     \\
b_0  & b_1 & \cdots & b_n  &  \cdots 
\end{array}
\right],   \  
{\rm then}  \  
(\alpha_0,\,\beta_0) \in 
D
\left(
\begin{array}{@{\,}cccc@{\,}}
a_0 & a_{1} & \cdots & a_n   \\
b_0 & b_{1} & \cdots & b_n  
\end{array}
\right)  \  ( n \geq 0). 
$
Conversely,  the Jacobi-Perron expansion of any point $(\alpha_0,\,\beta_0)$ $\in$   
$
D\!
\left(
\begin{array}{@{\,}cccc@{\,}}
a_0 & a_{1} & \cdots & a_n   \\
b_0 & b_{1} & \cdots & b_n  
\end{array}
\right)  \  ( n \geq 0) 
$\ 
has  the expansion from the first to $n+1$th step  
$ 
\left[
\begin{array}{@{\,}cccc@{\,}}
a_0  & a_1 & \cdots & a_n       \\
b_0  & b_1 & \cdots & b_n  
\end{array}
\right].  
$   

\noindent
Let $m (\geq 2)$ be a natural number. 
We consider a point 
 $(\alpha_0,\,\beta_0)$ $\in$ 
$
D\!
\left(
\begin{array}{@{\,}cccc@{\,}}
a_0 & a_{1} & \cdots & a_n   \\
b_0 & b_{1} & \cdots & b_n  
\end{array}
\right)   
$  
$( n \geq 0) $
satisfying  $\beta_{n+1} < m$. 
Since $\beta_{n+1} < m$ if and only if $\displaystyle{ 1  \over  m} < \{ \alpha_n \}$,  
we have 
$ 
(\{ \alpha_n \},\,\{ \beta_n \}) \in \mathcal{D}_{1/m} \ {\rm or}\ 
(\{ \alpha_n \},\,\{ \beta_n \}) \in \mathcal{D}^{\prime}_{1/m}
$  
according to $a_n < b_n$ or $a_n = b_n$. 
Therefore 
\[ 
\left\{ 
\left(
\begin{array}{@{\,}c@{\,}}
\alpha_0    \\
\beta_0   
\end{array} 
\right)  
\in 
D\!
\left(
\begin{array}{@{\,}cccc@{\,}}
a_0 & a_{1} & \cdots & a_n   \\
b_0 & b_{1} & \cdots & b_n  
\end{array}
\right)   
;\ \displaystyle{ 1  \over  m} < \{ \alpha_n \} 
\right\}  
\]  
\[ 
\subset 
\left\{ 
F
\left[
\begin{array}{@{\,}cccc@{\,}}
a_0  & a_1 & \cdots & a_{n}    \\
b_0  & b_1 & \cdots & b_{n}    
\end{array}
\right]  
\left(
\begin{array}{@{\,}c@{\,}}
X    \\
Y   
\end{array} 
\right)  
\ ;\ 
\left(
\begin{array}{@{\,}c@{\,}}
X    \\
Y   
\end{array} 
\right)  
\in 
\mathcal{D}_{1/m} \ \ {\rm or}\ \ \mathcal{D}^{\prime}_{1/m}  
\right\}. 
\] 
From the definitions, we have 
\[ 
\left|
\left\{ 
F
\left[
\begin{array}{@{\,}cccc@{\,}}
a_0  & a_1 & \cdots & a_{n}    \\
b_0  & b_1 & \cdots & b_{n}    
\end{array}
\right]  
\left(
\begin{array}{@{\,}c@{\,}}
X    \\
Y   
\end{array} 
\right)  
\ ;\ 
\left(
\begin{array}{@{\,}c@{\,}}
X    \\
Y   
\end{array} 
\right)  
\in 
\mathcal{D}_{1/m} \ \ {\rm or}\ \ \mathcal{D}^{\prime}_{1/m}  
\right\}  
\right|
\]  
\[ = 
S_n(1)-S_n(1/m)\ \ {\rm or}\ \ S_n^{\prime}(1)-S_n^{\prime}(1/m)
\ {\rm according\ to}\ a_n < b_n\ {\rm or}\ a_n = b_n. 
\] 

\bigskip

{\bf Theorem 4.2.}\ \  {\it Let        
$t\,  (0 < t < 1)$ be a real number and  $(a_n,\,b_n)\ (n \geq 0)$ a sequence with the admissibility, 
then we have }\ 
$ 
 S_n(t)   > t^2  S_n(1)  
\ \ {\it and}\ \ 
 S_n^{\prime}(t)     > t^2   S_n^{\prime}(1). 
$ 

\medskip

{\it Proof.}\ \  By means of Theorem 4.1, 
\[ 
\begin{array}{@{\,}l@{\,}}
\displaystyle{ S_n(t)   \over  S_n(1)  }
=
t^2 \cdot 
\displaystyle{ r_{n-2}+r_n   \over   r_{n-2}t+r_n  }
\cdot
\displaystyle{ r_{n-2}+r_{n-1}+r_n   \over   r_{n-2}t+r_{n-1}+r_n  }
\cdot 
\displaystyle{ r_{n-2}t+r_{n-1}+2r_n   \over   r_{n-2}t+r_{n-1}t+2r_n t  } > t^2,   \\  [13pt] 
\displaystyle{ S_n^{\prime}(t)   \over  S_n^{\prime}(1)  }
=
t \cdot 
\displaystyle{ r_{n-2}+r_{n-1}+r_n   \over   r_{n-2}t+r_{n-1}+r_n  }
\cdot 
\displaystyle{ r_{n-2}t+r_{n-1}+r_n (2-t)  \over   r_{n-2}t+r_{n-1}t+r_n  }  
> t  > t^2. \ \ \ \ \ \ \hspace{19mm}{\rm \qed}
\end{array}  
\]

\bigskip

Letting  $t=1/m$, we get the next result.   

\bigskip

{\bf Lemma 4.8.}\ \  {\it Let $m \geq 2$  be a natural number        
 and  $(a_n,\,b_n)\ (n \geq 0)$ a sequence with the admissibility,  
then we have }
\[ 
S_n(1)-S_n(1/m) < \left(1-\displaystyle{ 1  \over  m^2} \right) S_n(1)
\ \ {\it and}\ \ 
S_n^{\prime}(1)- S_n^{\prime}(1/m) < \left( 1-\displaystyle{ 1  \over  m^2} \right) S_n^{\prime}(1).
\] 

\bigskip

We define  for  integers $m \geq 2$ and $n \geq 0$, 
\[ 
\mathcal{D}_m(n) 
= \bigcup_{1 \leq b_i < m, \ 0 \leq a_i \leq b_i, \ 0 \leq i \leq n}
D\!
\left(
\begin{array}{@{\,}cccccc@{\,}}
a_0 & a_{1} & \cdots & a_i & \cdots & a_n   \\
b_0 & b_{1} & \cdots & b_i & \cdots & b_n  
\end{array}
\right)   
\] 
where we take the union of all  the $(a_i,\,b_i)$'s  satisfying the conditions. 

\bigskip

{\bf Theorem 4.3.}\ \  {\it It holds }       
\[ 
\left| \mathcal{D}_m(n+1) \right| < \left( 1-\displaystyle{ 1  \over  m^2} \right) \left| \mathcal{D}_m(n) \right|. 
\] 
 
\medskip

{\it Proof.}\ \  
\[  
 \mathcal{D}_m(n+1)  
= 
\displaystyle\bigcup_{1 \leq b_i < m, \ 0 \leq a_i \leq b_i, \ 0 \leq i \leq n+1}
D\!
\left(
\begin{array}{@{\,}cccc@{\,}}
a_0 & a_{1} & \cdots & a_{n+1}   \\
b_0 & b_{1} & \cdots & b_{n+1}  
\end{array}
\right)   
\] 

$ 
= 
\displaystyle\bigcup_{1 \leq b_i < m, \ 0 \leq a_i \leq b_i, \ 0 \leq i \leq n } 
\left( 
\displaystyle\bigcup_{1 \leq b_{n+1} < m, \ 0 \leq a_{n+1} \leq b_{n+1} } 
D\!
\left(
\begin{array}{@{\,}ccccc@{\,}}
a_0 & a_{1} & \cdots & a_n & a_{n+1}   \\
b_0 & b_{1} & \cdots & b_n & b_{n+1}  
\end{array}
\right) 
\right)  
$ 

\bigskip

$ 
 = 
\displaystyle\bigcup_{1 \leq b_i < m, \ 0 \leq a_i \leq b_i, \ 0 \leq i \leq n } 
\left\{ 
\left(
\begin{array}{@{\,}c@{\,}}
\alpha_0    \\
\beta_0   
\end{array} 
\right)  
\in 
D\!
\left(
\begin{array}{@{\,}cccc@{\,}}
a_0 & a_{1} & \cdots & a_n   \\
b_0 & b_{1} & \cdots & b_n  
\end{array}
\right)   
\ ;\ 
b_{n+1} < m 
\right\}  
$

\medskip
\noindent
(where $b_{n+1}=[ \beta_{n+1} ]$ is the second coordinate of the $n+2$th step of the expansion of $(\alpha_0,\,\beta_0)$)
\medskip 
 
$
= 
\displaystyle\bigcup_{1 \leq b_i < m, \ 0 \leq a_i \leq b_i, \ 0 \leq i \leq n } 
\left\{ 
\left(
\begin{array}{@{\,}c@{\,}}
\alpha_0    \\
\beta_0   
\end{array} 
\right)  
\in 
D\!
\left(
\begin{array}{@{\,}cccc@{\,}}
a_0 & a_{1} & \cdots & a_n   \\
b_0 & b_{1} & \cdots & b_n  
\end{array}
\right)   
\ ;\ 
\displaystyle{  1  \over  m } < \{ \alpha_n \}  
\right\}  
$ 

\medskip
\noindent
(here $\alpha_n$ comes from $\psi^n(\alpha_0,\,\beta_0)$)
\medskip

$
\subset 
 \displaystyle\bigcup_{1 \leq b_i < m, \ 0 \leq a_i \leq b_i, \ 0 \leq i \leq n } 
\left\{ 
F\!
\left[
\begin{array}{@{\,}cccc@{\,}}
a_0 & a_{1} & \cdots & a_n   \\
b_0 & b_{1} & \cdots & b_n  
\end{array}
\right]
\left(
\begin{array}{@{\,}c@{\,}}
X    \\
Y   
\end{array} 
\right)     
 ; 
\left(
\begin{array}{@{\,}c@{\,}}
X    \\
Y   
\end{array} 
\right)     
\in 
\mathcal{D}_{1/m} \ \ {\rm or}\ \ \mathcal{D}^{\prime}_{1/m}  
\right\}, 
$ 

\bigskip

\noindent 
according to $a_n <b_n$ or $a_n=b_n$, respectively. 
Therefore, we have  
\[  
\left| \mathcal{D}_m(n+1) \right| 
\leq 
\sum_{1 \leq b_i < m, \ 0 \leq a_i \leq b_i, \ 0 \leq i \leq n } 
\left\{
S_n(1)- S_n(1/m) 
\ \ \ {\rm or}\ \ \ 
S_n^{\prime}(1)- S_n^{\prime}(1/m) 
\right\}
\]   

$
< 
\displaystyle\sum_{1 \leq b_i < m, \ 0 \leq a_i \leq b_i, \ 0 \leq i \leq n } 
\left( 1-\displaystyle{ 1  \over  m^2} \right) 
\left\{ S_n(1) \ \ \ {\rm or}\ \ \ S_n^{\prime}(1) \right\}     
$ 

\bigskip

$
= 
\left( 1-\displaystyle{ 1  \over  m^2} \right) 
\displaystyle\sum_{1 \leq b_i < m, \ 0 \leq a_i \leq b_i, \ 0 \leq i \leq n } 
\left\{ S_n(1) \ \ \ {\rm or}\ \ \ S_n^{\prime}(1) \right\}  
$ 

\bigskip

$
= 
\left( 1-\displaystyle{ 1  \over  m^2} \right) 
\displaystyle\sum_{1 \leq b_i < m, \ 0 \leq a_i \leq b_i, \ 0 \leq i \leq n } 
\left|
D\!
\left(
\begin{array}{@{\,}cccc@{\,}}
a_0 & a_{1} & \cdots & a_n   \\
b_0 & b_{1} & \cdots & b_n  
\end{array}
\right)   
\right|  
$ 

\bigskip

$
= 
\left( 1-\displaystyle{ 1  \over  m^2} \right) 
\left|
 \mathcal{D}_m(n) 
\right|.  
$  
\bigskip

\noindent
Note that $S_n$ and $S_n^{\prime}$ depend on  the sequences $(a_i,\,b_i)$ $(0 \leq i \leq n)$  
satisfying the conditions about  summation and the admissibility.  
\hfill  \qed

\bigskip

We define  
\[ 
\begin{array}{@{\,}l@{\,}}
\mathcal{B}=\{ (\alpha_0, \beta_0) \in \mathcal{D}\  ; \ \beta_n ( n \geq 0)\ {\rm is\ bounded\ from\ above} \},   \\
\mathcal{B}_m=\{ (\alpha_0, \beta_0) \in \mathcal{D}\  ; \ \beta_n < m\ {\rm \ for\ all}\  n \geq 0  \}\ \ 
(m \geq 2 {\rm \ is\ a\ natural\ number}).  
\end{array}
\] 
Then,  
$ 
\mathcal{B} = \displaystyle\bigcup_{m \geq 2} \mathcal{B}_m, \  \  
\mathcal{B}_m \subset \bigcap_{n \geq 0}  \mathcal{D}_m(n),  
\  
{\rm and}\  \ 
\mathcal{B}_m \subset  \mathcal{D}_m(n)  \ \ {\rm for\ all}\ n \geq 0.
$ 

\bigskip

{\bf Theorem 4.4.}\ \ {\it  $\mathcal{B}$ is a null set.}        

\bigskip

{\it Proof.}\ \ 
Note that  
$ 
\left|
D\!
\left(
\begin{array}{@{\,}c@{\,}}
a_0    \\
b_0  
\end{array}
\right)   
\right| 
= 
1 \ {\rm or}\ 1/2,\ {\rm when}\ a_0<b_0\  {\rm or}\ a_0 = b_0,\ {\rm respectively}, 
$

\[ 
 \mathcal{D}_m(0) = \bigcup_{1 \leq b_0 < m,\ 0 \leq a_0 \leq b_0} 
D\!
\left(
\begin{array}{@{\,}c@{\,}}
a_0    \\
b_0  
\end{array}
\right)   \quad {\rm and} \quad \left| \mathcal{D}_m(0) \right| =(m^2-1)/2. 
\] 

By Theorem 4.3, ,  
\[ 
\left|
 \mathcal{D}_m(n) 
\right| 
< 
\left( 1-\displaystyle{ 1  \over  m^2} \right)^{\!n} 
\!\!\left|
 \mathcal{D}_m(0 )
 \right| 
 = 
\left( 1-\displaystyle{ 1  \over  m^2} \right)^{\!n} 
\!\displaystyle{ m^2-1  \over  2 }. 
\]  
Thus, 
$ 
\left|
 \mathcal{D}_m(n) 
\right| 
$\     
tends to $0$ when $n \longrightarrow \infty$. 
So, we have 
$\left| \mathcal{B}_m \right| =0$ for all $m \geq 2$ and therefore  $\left| \mathcal{B} \right| =0$. 
\hfill  \qed   

\bigskip 

This theorem means that, the set of all points ($\in  \mathcal{D}$) with periodic expansion, 
 a subset of $\mathcal{B}$, is null. 

In case of $m=2$, that is $\beta_n <2$ for all $n \geq 0$,  since $(a_n,\,b_n)$ can only take 
$(0,\,1)$ or $(1,\,1)$, it holds 
\[
\left( \!
\begin{array}{@{\,}cc@{\,}}
\alpha_0  \\
\beta_0  
\end{array}
\!\! \right) 
=
\left[
\begin{array}{@{\,}cccc@{\,}}
0 & 0 & 0 & \cdots    \\
1 & 1 & 1 & \cdots   
\end{array}
\right]  
\] 
or 
\[
\left( \!
\begin{array}{@{\,}cc@{\,}}
\alpha_0  \\
\beta_0  
\end{array}
\!\! \right) 
=
\left[
\begin{array}{@{\,}ccccccc@{\,}}
0 & 0 & \cdots & 0 & 1 & 1 & \cdots    \\
1 & 1 & \cdots & 1 & 1 & 1 & \cdots   
\end{array}
\right]  
\]
by the admissibility. 
Thus, $(\alpha_0,\,\beta_0)$ has the periodic expansion.

\vspace*{3\baselineskip}

\noindent
{\Large \bf References}

\medskip

\begin{description} 
\item[  [1\!\!] ]\  
Bernstein, L., 
{\it The Jacobi-Perron Algorithm---Its theory and Application}, 
Springer-Verlag, Berlin, Lecture Notes in Mathematics, \textbf{207} (1971). 

\item[ [2\!\!] ]\ 
Dubois, E., Farhane, A. and Paysant-Le Roux, R., 
{\it The Jacobi-Perron Algorithm and Pisot numbers}, 
 Acta Arith., \textbf{111} (2004), 269--275. 

\item[ [3\!\!] ]\ Hardy, G.H., and Wright, E.M., 
{\it An Introduction to the theory of numbers}, 
Oxford at the Clarendon Press, London (1983).  

\item[ [4\!\!] ]\ Jacobi, C.G.J., 
{\it Allgemeine Theorie der kettenbruchaehnlichen Algorithmen}, 
 J. reine angew. Math., \textbf{69} (1869), 29--64. 

\item[ [5\!\!] ]\ 
Nakaishi, K., 
{\it The Exponent of Convergence for 2-Dimensional Jacobi-Perron Type Algorithms}, 
 Monatsh. Math., \textbf{132} (2001), 141--152. 

\item[ [6\!\!] ]\ 
Paley, R.E.C.E. and Ursell, H.D., 
{\it Continued fractions in several dimensions}. 
 Proc. Cambridge Philos., \textbf{26} (1930), 127--144. 

\item[ [7\!\!] ]\ 
Perron, O., 
{\it Grundlagen f\"{u}r eine Theorie des Jacobischen Kettenbruchalgorithmus},  
 Math. Ann.,  \textbf{64} (1907), 1--76. 

\item[ [8\!\!] ]\ 
Schweiger, F., 
{\it Metrische S\"{a}tze \"{u}ber den Jacobischen Algorithmus}, 
 Monatsh. Math., \textbf{69} (1965), 243--255. 

\item[ [9\!\!] ]\ 
------, 
{\it Multidimensional Continued Fractions}, 
Oxford Univ. Press, Oxford, 2000.

\item[ [10\!\!] ]\ 
------, 
{\it Some remarks on Diophantine approximation by the Jacobi-Perron algorithm}, 
 Acta Arith., \textbf{133} (2008), 209--219.

\end{description} 

\vspace*{2\baselineskip}

\noindent
Tomizawa 4--17--40--601, Taihaku-ku, Sendai, 982-0032,  JAPAN

\noindent
e-mail : wjkmx643@ybb.ne.jp 

\end{document}